%% file: main.tex
\documentclass{amsart}

\usepackage{amsmath, amssymb, amscd, amsfonts, mathrsfs} 

\usepackage[T1]{fontenc} 
\usepackage[english]{babel}

\usepackage{tikz-cd} 
\usepackage{enumitem} 
\usepackage{adjustbox}

\usepackage[hide]{todo} 
\usepackage{csquotes} 

\usepackage[backend=biber,style=alphabetic]{biblatex} 
\usepackage[shortcuts]{extdash}

\usepackage{xcolor}
\usepackage[only,mapsfrom]{stmaryrd}

\usepackage[colorlinks=true]{hyperref} 
\usepackage{cleveref} 

\numberwithin{equation}{subsection}

\newtheorem{introthm}{Theorem}  

\newtheorem{thm}{Theorem}[subsection]
\newtheorem{prop}[thm]{Proposition}
\newtheorem{lem}[thm]{Lemma}
\newtheorem{cor}[thm]{Corollary}
\newtheorem{defi}[thm]{Definition}
\newtheorem{defiprop}[thm]{Definition-Proposition}

\newtheorem*{prop*}{Proposition}
\newtheorem*{thm*}{Theorem}
\newtheorem*{defi*}{Definition}

\theoremstyle{definition}
\newtheorem{ex}[thm]{Example}    

\theoremstyle{remark}
\newtheorem{rmk}[thm]{Remark}    

\crefname{thm}{thm.}{thms.}
\Crefname{thm}{Theorem}{Theorems}
\crefname{lem}{lem.}{lemmas}
\Crefname{lem}{Lemma}{Lemmas}
\crefname{prop}{prop.}{props.}
\Crefname{prop}{Proposition}{Propositions}
\crefname{defi}{def.}{defs.}
\Crefname{defi}{Definition}{Definitions}
\crefname{section}{§}{§}
\Crefname{section}{Section}{Sections}
\crefname{subsection}{§}{§}
\Crefname{subsection}{Subsection}{Subsections}
\crefname{paragraph}{§}{§}
\Crefname{paragraph}{Paragraph}{Paragraphs}

\newcommand*{\sheafhom}{\mathscr{H}\kern -1.5pt om}
\newcommand*{\scrLie}{\mathscr{L}\kern -1pt ie \kern 0.5pt}

\newcommand{\cal}{\mathcal}
\newcommand{\scr}{\mathscr}

\newcommand{\RR}{\mathbb{R}}
\newcommand{\ZZ}{\mathbb{Z}}

\newcommand{\QQ}{\mathbb{Q}}
\newcommand{\PP}{\mathbb{P}}
\newcommand{\GG}{\mathbb{G}}
\newcommand{\CC}{\mathbb{C}}
\DeclareMathOperator{\Spec}{Spec}
\DeclareMathOperator{\Hom}{Hom}

\DeclareMathOperator{\Ext}{Ext}

\DeclareMathOperator{\Sh}{Sh}

\DeclareMathOperator{\Lie}{Lie}
\DeclareMathOperator{\ord}{ord}

\DeclareMathOperator{\len}{length}
\DeclareMathOperator{\et}{et}

\mathchardef\mhyphen="2D 

\title{Constructible tori over Dedekind schemes}
\author{Adrien Morin and Takashi Suzuki}

\begin{document}


\begin{abstract}
	We introduce an exact category of torsion-free constructible tori and an abelian category of constructible tori over a Dedekind scheme with perfect residue fields. The first one has an explicit description as $2$-term complexes of smooth commutative group algebraic spaces. Using the second-named author's duality results \cite{Suz19}, we prove that they are equivalent to the opposite of the categories of torsion-free $\ZZ$-constructible sheaves and all $\ZZ$-constructible sheaves, respectively. We then define $L$-functions for constructible tori over a Dedekind scheme proper over $\ZZ$ in terms of their étale realizations and prove a special value formula at $s=0$ using the Weil-étale formalism developed by the first-named author in \cite{Morin2023b}. This extends the results of the first-named author by removing the tame ramification hypothesis.
\end{abstract}

\maketitle

\tableofcontents


\input{constTori}
\input{specialValue}

\printbibliography

\end{document}

%% file: constTori.tex
\subsection{Introduction}
Let $X$ be an irreducible Dedekind scheme with perfect residue fields. We define torsion-free constructible tori over $X$ as complexes
		\[
		0
		\to
		\mathcal{T}
		\to
		\bigoplus_{v \in X_{0}}
		i_{v, \ast} E_{v}
		\to
		0
		\]
		of commutative smooth group algebraic spaces over $X$
		concentrated in degrees $0$ and $1$,
		where the term $\mathcal{T}$ is the N\'eron (lft) model of a torus over $K$
		and the term $E_{v}$ is a torsion-free
		$\mathbb{Z}$-constructible sheaf over $v$,
		such that the induced morphism
		\begin{equation}
			\pi_{0}(\mathcal{T}_{v})
			\to
			E_{v}
		\end{equation}
		over $v$ is an isomorphism for almost all $v \in X_{0}$.

The category of torsion-free constructible tori has a natural realization functor to $D(X_{\mathrm{sm}})$, by considering the associated complex of representable sheaves on $X_\mathrm{sm}$ up to quasi-isomorphism. As an example, let $\mathcal{T}$ be the N\'eron model of a torus over $K$. The complex
	\[
	[\mathcal{T}
	\to
	\bigoplus_{v \in X_{0}}
	i_{v, \ast} \pi_{0}(\mathcal{T}_{v}) / \mathrm{tor}]
	\]
	is a torsion-free constructible torus over $X$, whose realization in $D(X_{\mathrm{sm}})$ identifies with the sheaf represented by the maximal open subgroup scheme $\mathcal{T}^{\tau}$ of $\mathcal{T}$ whose fibers at all closed points have torsion group of components.

One can associate a torsion-free constructible torus to a torsion-free $\ZZ$-con-structible sheaf: given a torsion-free $\ZZ$-constructible sheaf $F$ on $X$, there is a canonical map \[\mathcal{T}\to \bigoplus_{v \in X_{0}} i_{v, \ast} \sheafhom(F_v,\ZZ)\]
where $\mathcal{T}$ is the N\'eron model over $X$ of the Cartier dual of $F_{K}$, which we call the dual $F^D$ of $F$. Conversely, for $\mathscr{T}$ a torsion-free constructible torus over $X$, then 
		\[
		\mathscr{T}^{D}
		:=
		\operatorname{\mathbf{Ext}}^{0}_{X_{\mathrm{sm}}}(
		\mathscr{T},
		\mathbf{G}_{m}
		).
		\]
		is a torsion-free $\ZZ$-constructible sheaf on $X$.
		The dual of the constant sheaf $\ZZ$ is represented in $D(X_{\mathrm{sm}})$ by the split torus $\GG_m$. More generally, the dual of the extension by zero of a lattice $F$ on a dense open $U$ is (represented by) the Néron model of the Cartier dual $T_U$ of $F$ on $U$.

The image of $F^D$ in $D(X_{\mathrm{et}})$ identifies with the étale inner hom $R \operatorname{\mathbf{Hom}}_{X_{\mathrm{et}}}(F,\GG_m)$ (\cref{0026}). When $X$ is proper over $\ZZ$, Artin-Verdier duality says that the étale cohomology with appropriate compact support of a torsion-free $\ZZ$-constructible sheaf $F$ is dual to the étale cohomology of $F^D$, which is one justification of the name \emph{dual}. Moreover:

\begin{introthm}[\cref{0041}]
	Let $\mathbb{Z} \mhyphen \mathrm{Con}_{\mathrm{tf}} / X$ denote the category of torsion-free $\ZZ$-constructible sheaves and $\mathrm{CT}_{\mathrm{tf}} / X$ the category of torsion-free constructible tori. The functors
	\[\begin{tikzcd}[ampersand replacement=\&]
		{\mathbb{Z} \mhyphen \mathrm{Con}_{\mathrm{tf}} / X} \& {(\mathrm{CT}_{\mathrm{tf}} / X)^\mathrm{op}}
		\arrow["{(-)^D}", out= 10, in=170, from=1-1, to=1-2]
		\arrow["{(-)^D}", out=190, in=-10, from=1-2, to=1-1]
	\end{tikzcd}\]
	are quasi-inverse equivalences of categories.
\end{introthm}
This theorem induces a natural notion of quasi-isomorphism for $\mathrm{CT}_{\mathrm{tf}} / X$-valued cochain complexes, by asking that the cone has acyclic dual (as a complex of $\ZZ$-constructible sheaves). Forming the respective localizations, we get an equivalence of bounded derived categories
\[
D^b(\mathrm{CT}_{\mathrm{tf}} / X)\simeq
D^{b}(
\mathbb{Z} \mhyphen \mathrm{Con}_{\mathrm{tf}} / X
)^{\mathrm{op}}
\simeq
D^{b}(\mathbb{Z} \mhyphen \mathrm{Con} / X)^{\mathrm{op}}.
\]

We can thus identify an abelian category of constructible tori as the essential image of $(\mathbb{Z} \mhyphen \mathrm{Con} / X)^\mathrm{op}$ in $D^b(\mathrm{CT}_{\mathrm{tf}} / X)$ under the quasi-inverse of the above. More explicitly:
\begin{defi*}
	A \emph{constructible torus over $X$} is
	an object of $D^{b}(\mathrm{CT}_{\mathrm{tf}} / X)$
	represented by a two-term complex
	$0 \to \mathscr{T}' \to \mathscr{T}'' \to 0$
	in $\mathrm{CT}_{\mathrm{tf}} / X$
	concentrated in degrees $0$ and $1$ such that
	$\mathscr{T}_{K}' \to \mathscr{T}_{K}''$ is faithfully flat
	and $\mathscr{T}_{v}' \to \mathscr{T}_{v}''$ has
	finite cokernel for all $v \in X_{0}$.
	The constructible tori over $X$ form
	a full subcategory of $D^{b}(\mathrm{CT}_{\mathrm{tf}} / X)$,
	which we denote by $\mathrm{CT} / X$.
\end{defi*}

\begin{introthm}[\cref{0043}]
	The category of constructible tori over $X$ is contravariantly equivalent to the category of $\mathbb{Z}$-constructible sheaves over $X$.
\end{introthm}
The second-named author developed a notion of Cartier duality for Néron models of $1$-motives over $K$. In that setting, the Néron model of a lattice $F_K$ over $K$ is the complex $\tau^{\leq 1}Rg_\ast F_K\in D^b(\mathbb{Z} \mhyphen \mathrm{Con} / X)$, whose Cartier dual is the connected Néron model $\mathcal{T}^0$ of the Cartier dual of $F_K$. Thus we can interpret the object $(\tau^{\leq 1}Rg_\ast F_K)^D\in D^b(\mathrm{CT}/X)$ as being an incarnation of $\mathcal{T}^0$; in particular its étale realization is represented by $\mathcal{T}^0$. On the other hand, $(g_\ast F_K)^D$ is the torsion-free constructible torus $\mathcal{T}^\tau$ mentioned above.

%

The realization to $D(X_{\mathrm{sm}})$ does not pass to the bounded derived category of constructible tori, but the étale realization does, and one can thus speak of the étale cohomology of a bounded complex of constructible tori. In a second part, we apply those constructions to the study of the special value at $s=0$ of $L$-functions associated \emph{via} their étale realization to constructible tori over a Dedekind scheme proper over $\Spec(\ZZ)$:
\begin{defi*}
	Let $\scr{T}^\bullet \in D^b(\mathrm{CT}/X)$. We define the $L$-function of $\scr{T}^\bullet$ by the Euler product
	\[
	L_X(\scr{T}^\bullet,s)=\prod_{x\in X_0} \det\bigl (I-N(x)^{-s}\varphi|\bigl (i_x^\ast \scr{T}^\bullet_{\mathrm{et}}\bigr )\widehat{\otimes}\QQ_{\ell_x}\bigr )^{-1}
	\]
	where $\ell_x$ is coprime to the residual characteristic at $x$ and $\varphi$ is the \emph{geometric} Frobenius at $x$.
\end{defi*}
Those $L$-functions are essentially Artin $L$-functions of rational representations in the variable $s+1$, and we prove a special value formula generalizing that of the first-named author in \cite{Morin2023b}. In particular, we remove the tame ramification hypothesis from \cite{Morin2023b}. The improvement follows from the results in the first part together with a multiplicativity formula for the determinant of the Lie algebra of the Néron model of a torus (\cref{mult_det_lie_OX}), which implies that the functor associating to a torus the determinant of the Lie algebra of its Néron model induces a functor on constructible tori:
\begin{prop*}[\cref{prop:extension_detLie}]
	The assignment $\mathscr{T} \mapsto \det_{\ZZ} \operatorname{Lie}(\mathscr{T})$, valued in graded $\ZZ$-lines, extends uniquely to a functor $\Delta_X^{\mathrm{add}}:D^b(\mathrm{CT}/X)_{\mathrm{iso}} \to \cal{P}_{\ZZ}$ from the groupoid core of $D^b(\mathrm{CT}/X)$ to the Picard groupoid of graded $\ZZ$-lines.
\end{prop*}
We use this to define a Weil-étale Euler characteristic $\chi_X$ on $D^b(\mathrm{CT}/X)$ incorporating the correct ``additive factor'', following the approach in \cite{Morin2023b}. The special value formula is then:
\begin{introthm}[\cref{thm:special_value}]
	Let $X$ be a Dedekind scheme proper over $\Spec(\ZZ)$ and let $\scr{T}^\bullet\in D^b(\mathrm{CT}/X)$. We have the special value formula
	\[
	L_X^\ast(\scr{T}^\bullet,0)=\pm \chi_X(\scr{T}^\bullet).
	\]
\end{introthm}
An unconditional formula for the vanishing order was already given in \cite{Morin2023b}. As an application, we obtain a special value formula at $s=1$ for Artin $L$-functions of rational representations, as explained in \cite{Morin2023b}.

\subsubsection{Acknowledgments}
The first-named author would like to thank Baptiste Morin and Fabien Pazuki for their encouragements and support.

\subsection{Torsion-free constructible tori}

We will define ``torsion-free constructible tori'' generalizing tori
in a way parallel to how we define torsion-free
$\mathbb{Z}$-constructible sheaves generalizing lattices.

Let $X$ be a Dedekind scheme (a noetherian regular scheme of dimension $\le 1$).
Assume that $X$ is irreducible and
the residue field $\kappa(v)$ of $X$ at any closed point $v$ is perfect.
(The latter assumption is needed to use the duality results in \cite{Suz19}.)
Let $X_{0}$ be the set of closed points of $X$.
For each $v \in X_{0}$, let $i_{v} \colon v \hookrightarrow X$ be the inclusion.
Let $K$ be the function field of $X$.
For an $X$-scheme $Y$,
set $Y_{v} = Y \times_{X} v$ for each $v \in X_{0}$
and $Y_{K} = Y \times_{X} \operatorname{Spec} K$.
For a commutative group scheme $G$ locally of finite type over $v$,
let $G^{0} \subset G$ be the identity component
and $\pi_{0}(G) = G / G^{0}$ the \'etale group scheme
of components (over $v$).
Now we define our main objects:

\begin{defi}
	A \emph{torsion-free constructible torus over $X$} is a complex
		\[
				0
			\to
				\mathcal{T}
			\to
				\bigoplus_{v \in X_{0}}
					i_{v, \ast} E_{v}
			\to
				0
		\]
	of commutative smooth group algebraic spaces over $X$
	concentrated in degrees $0$ and $1$,
	where the term $\mathcal{T}$ is the N\'eron (lft) model of a torus over $K$
	and the term $E_{v}$ is a torsion-free
	$\mathbb{Z}$-constructible sheaf over $v$,
	such that the induced morphism
		\begin{equation} \label{0000}
				\pi_{0}(\mathcal{T}_{v})
			\to
				E_{v}
		\end{equation}
	over $v$ is an isomorphism for almost all $v \in X_{0}$.
	
	A morphism of torsion-free constructible tori is a morphism of complexes
	of group algebraic spaces.
\end{defi}

We denote a torsion-free constructible torus using the mapping cone notation
$[\mathcal{T} \to \bigoplus_{v \in X_{0}} i_{v, \ast} E_{v}][-1]$
(where $[-1]$ is the shift, so $[\;\cdot\;][-1]$ is the mapping fiber).

Torsion-free constructible tori form an additive category
but do not form an abelian category.
It is occasionally useful to work
within a category of sheaves containing them.
Let $X_{\mathrm{sm}}$ be the category of smooth $X$-schemes
with $X$-scheme morphisms endowed with the \'etale topology.
Let $\operatorname{Ab}(X_{\mathrm{sm}})$ be
the category of sheaves of abelian groups over $X_{\mathrm{sm}}$
with category of complexes $\mathrm{Ch}(X_{\mathrm{sm}})$
and derived category $D(X_{\mathrm{sm}})$.
Let $\operatorname{\mathbf{Hom}}_{X_{\mathrm{sm}}}$
and $\operatorname{\mathbf{Ext}}^{n}_{X_{\mathrm{sm}}}$ for each $n$
be the sheaf-Hom and $n$-th sheaf-Ext functors, respectively,
for $\operatorname{Ab}(X_{\mathrm{sm}})$.
An \'etale morphism $f \colon Y \to X$ from another (Dedekind) scheme
induces a morphism of sites
$f \colon Y_{\mathrm{sm}} \to X_{\mathrm{sm}}$
defined by the functor $X' \times_{X} Y \mapsfrom X'$
on the underlying categories.

The category of torsion-free constructible tori can be viewed
as a full subcategory of $\mathrm{Ch}(X_{\mathrm{sm}})$.
Hence a torsion-free constructible torus defines an object of
$\mathrm{Ch}(X_{\mathrm{sm}})$ and hence of $D(X_{\mathrm{sm}})$. It will follow from the biduality established later that the functor from torsion-free constructible tori to $D(X_{\mathrm{sm}})$ is faithful.

\begin{ex} \label{0045}
	Let $\mathcal{T}$ be the N\'eron model of a torus over $K$.
	For each $v \in X_{0}$, let $\pi_{0}(\mathcal{T}_{v}) / \mathrm{tor}$ be
	the torsion-free quotient of $\pi_{0}(\mathcal{T}_{v})$.
	Then the canonical morphism
		\[
				\mathscr{T}
			\colon
				\mathcal{T}
			\to
				\bigoplus_{v \in X_{0}}
					i_{v, \ast} \pi_{0}(\mathcal{T}_{v}) / \mathrm{tor}
		\]
	viewed as a complex in degrees $0$ and $1$
	is a torsion-free constructible torus over $X$.
	Note that this morphism is smooth and faithfully flat.
	Its kernel is the relative torsion component $\mathcal{T}^{\tau}$
	(see also \cite[Section 6]{Kle05}),
	namely the maximal open subgroup scheme of $\mathcal{T}$
	with torsion $\pi_{0}(\mathcal{T}^{\tau}_{v})$ for all $v \in X_{0}$.
	Hence the image of $\mathscr{T}$ in $D(X_{\mathrm{sm}})$
	can be identified with $\mathcal{T}^{\tau}$.
\end{ex}

We define a notion of good reduction of torsion-free constructible tori.
Note that for a torus $\mathcal{T}_{U}^{0}$
over a dense open subscheme $j_{U} \colon U \hookrightarrow X$
with generic fiber $T$,
we have two canonical objects over $X$:
the N\'eron model $\mathcal{T}$ over $X$ of $T$ over $K$,
and the N\'eron model $\mathcal{T}^{(U)}$ over $X$
of $\mathcal{T}_{U}^{0}$ over $U$.
The latter represents $j_{U, \ast} \mathcal{T}_{U}^{0}$
and is the open subgroup scheme of the former with connected fibers over $U$.
Hence we have an exact sequence
	\begin{equation} \label{0022}
			0
		\to
			\mathcal{T}^{(U)}
		\to
			\mathcal{T}
		\to
			\bigoplus_{v \in U_{0}}
				i_{v, \ast} \pi_{0}(\mathcal{T}_{v})
		\to
			0.
	\end{equation}

\begin{prop} \label{0024}
	Let $\mathscr{T} = [\mathcal{T} \to \bigoplus_{v \in X_{0}} i_{v, \ast} E_{v}][-1]$
	be a torsion-free constructible torus over $X$.
	Let $U \subset X$ be a dense open subscheme such that
	the relative identity component $\mathcal{T}_{U}^{0}$ (\cite[\href{https://stacks.math.columbia.edu/tag/055R}{Tag 055R}]{Sta22})
	of $\mathcal{T}_{U}$ is a torus.
	Let $\mathcal{T}^{(U)}$ be the N\'eron model over $X$ of $\mathcal{T}_{U}^{0}$.
	Then we have a commutative diagram
		\[
			\begin{CD}
					0
				@>>>
					\mathcal{T}^{(U)}
				@>>>
					\mathcal{T}
				@>>>
					\bigoplus_{v \in U_{0}}
						\pi_{0}(\mathcal{T}_{v})
				@>>>
					0
				\\
				@.
				@VVV
				@VVV
				@VVV
				@.
				\\
					0
				@>>>
					\bigoplus_{v \in X \setminus U}
						i_{v, \ast} E_{v}
				@>>>
					\bigoplus_{v \in X_{0}}
						i_{v, \ast} E_{v}.
				@>>>
					\bigoplus_{v \in U_{0}}
						i_{v, \ast} E_{v}
				@>>>
					0
			\end{CD}
		\]
	with exact rows.
\end{prop}

\begin{proof}
	This follows from \eqref{0022}.
\end{proof}

In particular, if moreover \eqref{0000} is an isomorphism
for all $v \in U_{0}$,
then we have a quasi-isomorphism
	\[
			\left[
					\mathcal{T}^{(U)}
				\to
					\bigoplus_{v \in X \setminus U}
						i_{v, \ast} E_{v}
			\right][-1]
		\to
			\left[
					\mathcal{T}
				\to
					\bigoplus_{v \in X_{0}}
						i_{v, \ast} E_{v}
			\right][-1]
	\]
in $\operatorname{Ch}(X_{\mathrm{sm}})$.

\begin{defi}
	Let
		$
				\mathscr{T}
			=
				[
						\mathcal{T}
					\to
						\bigoplus_{v \in X_{0}}
							i_{v, \ast} E_{v}
				][-1]
		$
	be a torsion-free constructible torus over $X$.
	Let $U \subset X$ be a dense open subscheme.
	\begin{enumerate}
		\item
			If the relative identity component
			$\mathcal{T}_{U}^{0}$ of $\mathcal{T}_{U}$ is a torus
			and \eqref{0000} is an isomorphism for all $v \in U_{0}$,
			we say that $\mathscr{T}$ \emph{has good reduction over $U$}.
		\item
			In this case, we call the complex
				\[
						\mathscr{T}^{(U)}
					=
						\left[
								\mathcal{T}^{(U)}
							\to
								\bigoplus_{v \in X \setminus U}
									i_{v, \ast} E_{v}
						\right][-1]
				\]
			the \emph{$U$-model of $\mathscr{T}$}.\footnote{To clarify, under our definitions this is \emph{not} a torsion-free constructible torus, but $\mathscr{T}$ and $\mathscr{T}^{(U)}$ are isomorphic in $D(X_{\mathrm{sm}})$.}
	\end{enumerate}
\end{defi}


\subsection{%
	\texorpdfstring{Duals of torsion-free $\mathbb{Z}$-constructible sheaves}
	{Duals of torsion-free Z-constructible sheaves}
}

We will construct a torsion-free constructible torus
contravariantly and canonically
from a torsion-free $\mathbb{Z}$-con-structible sheaf.
It will be defined as a certain combination
of the Cartier dual of a lattice over $U$
and the $\mathbb{Z}$-duals of lattices over $v \in X \setminus U$.
For this, we first need to calculate some Hom and Ext sheaves of lattices:

\begin{prop} \label{0006} \mbox{}
	\begin{enumerate}
		\item \label{0004}
			Let $j_{U} \colon U \hookrightarrow X$ be a dense open subscheme.
			Let $F_{U}$ be a lattice over $U$
			(a torsion-free locally constant $\mathbb{Z}$-constructible sheaf).
			Then the inner hom
				$				\operatorname{\mathbf{Hom}}_{X_{\mathrm{sm}}}(
						j_{U, !} F_{U},
						\mathbf{G}_{m}
					)
				$
			is represented by
			the N\'eron model over $X$ of the Cartier dual of $F_{U}$.
		\item \label{0005}
			Let $v \in X_{0}$.
			Let $F_{v}$ be a torsion-free $\mathbb{Z}$-constructible sheaf over $v$.
			Let $E_{v}$ be its $\mathbb{Z}$-dual.
			Then
				\begin{gather*}
							\operatorname{\mathbf{Hom}}_{X_{\mathrm{sm}}}(
								i_{v, \ast} F_{v},
								\mathbf{G}_{m}
							)
						=
							0,
					\\
							\operatorname{\mathbf{Ext}}^{1}_{X_{\mathrm{sm}}}(
								i_{v, \ast} F_{v},
								\mathbf{G}_{m}
							)
						\cong
							i_{v, \ast} E_{v}.
				\end{gather*}
	\end{enumerate}
\end{prop}

\begin{proof}
	\eqref{0004}
	This follows from
		\[
				\operatorname{\mathbf{Hom}}_{X_{\mathrm{sm}}}(
					j_{U, !} F_{U},
					\mathbf{G}_{m}
				)
			\cong
				j_{U, \ast}
				\operatorname{\mathbf{Hom}}_{U_{\mathrm{sm}}}(
					F_{U},
					\mathbf{G}_{m}
				)
		\]
	and the fact that
		$
				j_{U, \ast}
			\colon
				\operatorname{Ab}(U_{\mathrm{sm}})
			\to
				\operatorname{Ab}(X_{\mathrm{sm}})
		$
	is the N\'eron model functor by definition.
	
	\eqref{0005}
	Since
		$
				E_{v}
			\cong
				\operatorname{\mathbf{Hom}}_{v_{\mathrm{sm}}}(
					F_{v},
					\mathbb{Z}
				)
		$,
	we have
		$
				i_{v, \ast} E_{v}
			\cong
				\operatorname{\mathbf{Hom}}_{X_{\mathrm{sm}}}(
					i_{v, \ast} F_{v},
					i_{v, \ast} \mathbb{Z}
				)
		$.
	Let $U = X \setminus \{v\}$
	with the inclusion map $j_{U} \colon U \hookrightarrow X$.
	We have an exact sequence
		\begin{equation} \label{0010}
				0
			\to
				\mathbf{G}_{m}
			\to
				j_{U, \ast} \mathbf{G}_{m}
			\to
				i_{v, \ast} \mathbb{Z}
			\to
				0
		\end{equation}
	of commutative smooth group schemes over $X$.
	This induces a morphism
		\[
				\operatorname{\mathbf{Hom}}_{X_{\mathrm{sm}}}(
					i_{v, \ast} F_{v},
					i_{v, \ast} \mathbb{Z}
				)
			\to
				\operatorname{\mathbf{Ext}}^{1}_{X_{\mathrm{sm}}}(
					i_{v, \ast} F_{v},
					\mathbf{G}_{m}
				).
		\]
	To show it is an isomorphism,
	it is enough to show that
		$
				\operatorname{\mathbf{Ext}}^{n}_{X_{\mathrm{sm}}}(
					i_{v, \ast} F_{v},
					j_{U, \ast} \mathbf{G}_{m}
				)
			=
				0
		$
	for $n = 0$ and $1$.
	But this follows from
		\[
				R \operatorname{\mathbf{Hom}}_{X_{\mathrm{sm}}}(
					i_{v, \ast} F_{v},
					R j_{U, \ast} \mathbf{G}_{m}
				)
			\cong
				R j_{U, \ast}
				R \operatorname{\mathbf{Hom}}_{U_{\mathrm{sm}}}(
					j_{U}^{\ast} i_{v, \ast} F_{v},
					\mathbf{G}_{m}
				)
			=
				0.
		\]
\end{proof}

Now we will define the dual of a torsion-free $\mathbb{Z}$-constructible sheaf.
Let $F$ be a torsion-free $\mathbb{Z}$-constructible sheaf over $X$.
Let $\mathcal{T}$ be the N\'eron model over $X$ of the Cartier dual of $F_{K}$.
For each $v \in X_{0}$, let $E_{v}$ be the $\mathbb{Z}$-dual of $F_{v}$.

Let $j_{U} \colon U \hookrightarrow X$ be
a dense open subscheme where $F$ is locally constant.
Let $\mathcal{T}^{(U)}$ is the N\'eron model of the Cartier dual of $F_{U}$.
Consider the exact sequence
	\begin{equation} \label{0009}
			0
		\to
			j_{U, !} F_{U}
		\to
			F
		\to
			\bigoplus_{v \in X \setminus U}
				i_{v, \ast} F_{v}
		\to
			0
	\end{equation}
and the induced morphism
	\begin{equation} \label{0032}
			\operatorname{\mathbf{Hom}}_{X_{\mathrm{sm}}}(
				j_{U, !} F_{U},
				\mathbf{G}_{m}
			)
		\to
			\bigoplus_{v \in X \setminus U}
				\operatorname{\mathbf{Ext}}^{1}_{X_{\mathrm{sm}}}(
					i_{v, \ast} F_{v},
					\mathbf{G}_{m}
				).
	\end{equation}
Note by Proposition \ref{0006} \eqref{0005}
that this is part of the long exact sequence
	\begin{equation} \label{0036}
	\adjustbox{center}{\resizebox{\textwidth}{!}{$
	\begin{CD}
			@.
				0
			@>>>
				\operatorname{\mathbf{Hom}}_{X_{\mathrm{sm}}}(
					F,
					\mathbf{G}_{m}
				)
			@>>>
				\operatorname{\mathbf{Hom}}_{X_{\mathrm{sm}}}(
					j_{U, !} F_{U},
					\mathbf{G}_{m}
				)
			\\
			@>>>
				\displaystyle
				\bigoplus_{v \in X \setminus U}
					\operatorname{\mathbf{Ext}}^{1}_{X_{\mathrm{sm}}}(
						i_{v, \ast} F_{v},
						\mathbf{G}_{m}
					)
			@>>>
				\operatorname{\mathbf{Ext}}^{1}_{X_{\mathrm{sm}}}(
					F,
					\mathbf{G}_{m}
				)
			@>>>
				\operatorname{\mathbf{Ext}}^{1}_{X_{\mathrm{sm}}}(
					j_{U, !} F_{U},
					\mathbf{G}_{m}
				).
		\end{CD}$}}
	\end{equation}
Write \eqref{0032} as
	\begin{equation} \label{0008}
			\mathcal{T}^{(U)}
		\to
			\bigoplus_{v \in X \setminus U}
				i_{v, \ast} E_{v}
	\end{equation}
using Proposition \ref{0006}.

\begin{prop} \label{0012}
	In the above situation,
	if $F$ is constant over an \'etale neighborhood
	of a point $v \in X \setminus U$,
	then the morphism \eqref{0008} induces an isomorphism
	$\pi_{0}(\mathcal{T}_{v}) \overset{\sim}{\to} E_{v}$.
\end{prop}

\begin{proof}
	We may replace $X$ by its strict henselian localization
	$\operatorname{Spec} \mathcal{O}_{v}^{\mathrm{sh}}$ at $v$
	(so $U = \operatorname{Spec} K_{v}^{\mathrm{sh}}$).
	Hence we may assume that $F = \mathbb{Z}$.
	Then the first few terms of the long exact sequence for
		$
			\operatorname{\mathbf{Ext}}^{\;\cdot\;}_{X_{\mathrm{sm}}}(
				\;\cdot\;,
				\mathbf{G}_{m}
			)
		$
	associated with \eqref{0009} can be identified with \eqref{0010}.
	The induced morphism
	$\pi_{0}((j_{U, \ast} \mathbf{G}_{m})_{v}) \to \mathbb{Z}$
	is an isomorphism,
	which implies the result.
\end{proof}

For another dense open subscheme
$j_{V} \colon V \hookrightarrow X$ with $V \subset U$,
we have a commutative diagram
	\[
		\begin{CD}
				0
			@>>>
				j_{U, !} F_{U}
			@>>>
				F
			@>>>
				\bigoplus_{v \in X \setminus U}
					i_{v, \ast} F_{v}
			@>>>
				0
			\\ @. @V \mathrm{incl} VV @| @VV \mathrm{proj} V \\
				0
			@>>>
				j_{V, !} F_{V}
			@>>>
				F
			@>>>
				\bigoplus_{v \in X \setminus V}
					i_{v, \ast} F_{v}
			@>>>
				0.
		\end{CD}
	\]
This induces a commutative diagram
	\[
		\begin{CD}
				\mathcal{T}^{(U)}
			@>>>
				\bigoplus_{v \in X \setminus U}
					i_{v, \ast} E_{v}
			\\ @V \mathrm{incl} VV @VV \mathrm{incl} V \\
				\mathcal{T}^{(V)}
			@>>>
				\bigoplus_{v \in X \setminus V}
					i_{v, \ast} E_{v}.
		\end{CD}
	\]

\begin{prop} \label{0023}
	The above diagram fits in the commutative diagram
		\[
			\begin{CD}
					0
				@>>>
					\mathcal{T}^{(U)}
				@>>>
					\mathcal{T}^{(V)}
				@>>>
					\bigoplus_{v \in U \setminus V}
						\pi_{0}(\mathcal{T}_{v})
				@>>>
					0
				\\
				@.
				@VVV
				@VVV
				@VV \wr V
				@.
				\\
					0
				@>>>
					\bigoplus_{v \in X \setminus U}
						i_{v, \ast} E_{v}
				@>>>
					\bigoplus_{v \in X \setminus V}
						i_{v, \ast} E_{v}.
				@>>>
					\bigoplus_{v \in U \setminus V}
						i_{v, \ast} E_{v}
				@>>>
					0
			\end{CD}
		\]
	with exact rows.
\end{prop}

\begin{proof}
	This follows from \eqref{0022} and Proposition \ref{0012}.
\end{proof}

By the above manner, \eqref{0008} forms a filtered direct system in $U$.
Its direct limit defines a morphism
	\begin{equation} \label{0011}
			\mathcal{T}
		\to
			\bigoplus_{v \in X_{0}}
				i_{v, \ast} E_{v}
	\end{equation}
of group algebraic spaces over $X$.

\begin{prop}
	We have a commutative diagram
		\[
			\begin{CD}
					0
				@>>>
					\mathcal{T}^{(U)}
				@>>>
					\mathcal{T}
				@>>>
					\bigoplus_{v \in U_{0}}
						\pi_{0}(\mathcal{T}_{v})
				@>>>
					0
				\\
				@.
				@VVV
				@VVV
				@VV \wr V
				@.
				\\
					0
				@>>>
					\bigoplus_{v \in X \setminus U}
						i_{v, \ast} E_{v}
				@>>>
					\bigoplus_{v \in X_{0}}
						i_{v, \ast} E_{v}.
				@>>>
					\bigoplus_{v \in U_{0}}
						i_{v, \ast} E_{v}
				@>>>
					0
			\end{CD}
		\]
	with exact rows.
\end{prop}

\begin{proof}
	Take the direct limit in $V$ of the diagram in Proposition \ref{0023}.
\end{proof}

Hence \eqref{0008} and \eqref{0011} are quasi-isomorphic
as complexes in $\operatorname{Ab}(X_{\mathrm{sm}})$.
Here is the definition of the dual of $F$:

\begin{defi}
	Let $F$ be a torsion-free $\mathbb{Z}$-constructible sheaf over $X$.
	Define $F^{D}$ to be the morphism \eqref{0011}
	viewed as a complex concentrated in degrees $0$ and $1$.
	
	For a morphism $F \to G$
	of torsion-free $\mathbb{Z}$-constructible sheaves over $X$,
	define a morphism $G^{D} \to F^{D}$ of complexes of group algebraic spaces
	by the naturally induced morphism of the functoriality of \eqref{0009}
	(with $U$ chosen so that both $F$ and $G$ are locally constant).
\end{defi}

\begin{prop}
	Let $F$ be a torsion-free $\mathbb{Z}$-constructible sheaf over $X$.
	Then $F^{D}$ is a torsion-free constructible torus over $X$.
	If $F$ is locally constant over a dense open subscheme $U$,
	then $F^{D}$ has good reduction over $U$.
\end{prop}

\begin{proof}
	We only need to show that the induced morphism
	$\pi_{0}(\mathcal{T}_{v}) \to E_{v}$
	is an isomorphism for almost all $v \in X_{0}$.
	But this follows from Proposition \ref{0012}.
\end{proof}

Thus we have a contravariant functor $F \mapsto F^{D}$
from the category of torsion-free $\mathbb{Z}$-constructible sheaves over $X$
to the category of torsion-free constructible tori over $X$.
It is an additive functor.
We will give a non-trivial example of a dual in Example \ref{0048}.


\subsection{%
	\texorpdfstring{Comparison with derived Hom to $\mathbf{G}_{m}$}
	{Comparison with derived Hom to Gm}
}

We will clarify the relation between our dual $F^{D}$, the $\mathbf{G}_{m}$-dual $R \operatorname{\mathbf{Hom}}_{X_{\mathrm{sm}}}(F, \mathbf{G}_{m})$
 of $F$ on $X_{\mathrm{sm}}$ and its restriction
$R \operatorname{\mathbf{Hom}}_{X_{\mathrm{et}}}(F, \mathbf{G}_{m})$
to the small \'etale site of $X$.

Let $f \colon X_{\mathrm{sm}} \to X_{\mathrm{et}}$ be the morphism of sites
defined by the inclusion functor on the underlying categories.
Its pushforward functor $f_{\ast}$ is exact.
It sends a sheaf $F$ on $X_{\mathrm{sm}}$
to its restriction to $X_{\mathrm{et}}$.
For an object $F \in D(X_{\mathrm{sm}})$
(or a morphism between such objects),
we call $f_{\ast} F$ the restriction of $F$ to $X_{\mathrm{et}}$.
Let
$g \colon \operatorname{Spec} K \hookrightarrow X$
be the inclusion,
inducing a morphism of sites
$g \colon \operatorname{Spec} K_{\mathrm{sm}} \to X_{\mathrm{sm}}$
on the smooth sites (\cite[Proposition 2.1]{Suz19}).

The sheaf
$R^{1} g_{\ast} T \in \operatorname{Ab}(X_{\mathrm{sm}})$
is non-zero in general
as noted in \cite[Remark (4.8) (b)]{Cha00}.
See also \cite[Remarks (4.4) and (4.6)]{Cha00}.
The point of $F^{D}$ is to discard this uncontrollable sheaf
and cut out the representable parts:

\begin{prop} \label{0026}
	Let $F$ be a torsion-free $\mathbb{Z}$-constructible sheaf over $X$.
	Let $T$ be the Cartier dual of $F_{K}$.
	\begin{enumerate}
		\item \label{0020}
			There exists a canonical morphism
				\begin{equation} \label{0017}
						F^{D}
					\to
						R \operatorname{\mathbf{Hom}}_{X_{\mathrm{sm}}}(
							F,
							\mathbf{G}_{m}
						)
				\end{equation}
			in $D(X_{\mathrm{sm}})$.
			It is an isomorphism in $H^{0}$
			and fits in an exact sequence
				\begin{equation} \label{0019}
						0
					\to
						H^{1} F^{D}
					\to
						\operatorname{\mathbf{Ext}}^{1}_{X_{\mathrm{sm}}}(
							F,
							\mathbf{G}_{m}
						)
					\to
						R^{1} g_{\ast} T
				\end{equation}
			in $H^{1}$.
		\item \label{0021}
			The restriction of \eqref{0017} to $X_{\mathrm{et}}$ is an isomorphism.
	\end{enumerate}
\end{prop}

This almost follows from \eqref{0036}
modulo some technicalities about the choice of homotopy and mapping cones.
\footnote{
	Those familiar with $\infty$-categories	may avoid these technicalities. Namely, choosing an $U$-model of $F^D$ as in \ref{0008}, we have in the $\infty$-category $D(X_{\mathrm{sm}})$ an identification $F^D\simeq \mathrm{cofiber}(\tau_{\leq 1} (R\operatorname{\mathbf{Hom}}_{X_{\mathrm{sm}}}(j_{U,!}F_U,\GG_m)[-1]) \to \tau_{\leq 1} R\operatorname{\mathbf{Hom}}_{X_{\mathrm{sm}}}(\bigoplus_{v\in X\backslash U} i_{v,\ast} F_v,\GG_m))$. The canonical map is then the map induced on cofibers by the naturality of $\tau_{\leq 1}$. All of this is compatible with passing to colimits on $U$ by commutation of $\tau_{\leq 1}$ and cofibers with filtered colimits.
}

\begin{proof}
	\eqref{0020}
	Let $U \subset X$ range over
	dense open subschemes where $F$ is locally constant.
	By Proposition \ref{0006} \eqref{0005},
	we can apply the general Lemma \ref{0018} \eqref{0015} below
	to the bifunctor
		$
				\operatorname{\mathbf{Hom}}_{X_{\mathrm{sm}}}
			\colon
					\operatorname{Ab}(X_{\mathrm{sm}})^{\mathrm{op}}
				\times
					\operatorname{Ab}(X_{\mathrm{sm}})
			\to
				\operatorname{Ab}(X_{\mathrm{sm}})
		$
	(where $\mathrm{op}$ denotes the opposite category),
	the object $\mathbf{G}_{m}$ and the exact sequence \eqref{0009} for varying $U$.
	This produces the morphism \eqref{0017}.
	By Lemma \ref{0018} \eqref{0014} below,
	it is an isomorphism in $H^{0}$ and
	the direct limit in $U$ of the monomorphisms
		\[
				\operatorname{Coker} \left(
						\operatorname{\mathbf{Hom}}_{X_{\mathrm{sm}}}(
							j_{U, !} F_{U},
							\mathbf{G}_{m}
						)
					\to
						\bigoplus_{v \in X \setminus U}
							\operatorname{\mathbf{Ext}}^{1}_{X_{\mathrm{sm}}}(
								i_{v, \ast} F_{v},
								\mathbf{G}_{m}
							)
				\right)
			\hookrightarrow
				\operatorname{\mathbf{Ext}}^{1}_{X_{\mathrm{sm}}}(
					F,
					\mathbf{G}_{m}
				)
		\]
	coming from \eqref{0009}.
	The cokernel of this monomorphism for each $U$ injects into the sheaf
		\[
				\operatorname{\mathbf{Ext}}^{1}_{X_{\mathrm{sm}}}(
					j_{U, !} F_{U},
					\mathbf{G}_{m}
				)
			\cong
				R^{1} j_{U, \ast} F_{U}.
		\]
	Its direct limit in $U$ is $R^{1} g_{\ast} T$,
	giving \eqref{0019}.
	
	\eqref{0021}
	From the proof of \eqref{0020}, it is enough to show that
	the restrictions of the sheaves
	$R^{n} g_{\ast} T$ and
		$
			\operatorname{\mathbf{Ext}}^{n + 1}_{X_{\mathrm{sm}}}(
				i_{v, \ast} F_{v},
				\mathbf{G}_{m}
			)
		$
	to $X_{\mathrm{et}}$ for any $n \ge 1$ and $v \in X_{0}$ are zero.
	The stalk of $R^{n} g_{\ast} T$ at $v$ as an \'etale sheaf is
	$H^{n}(K_{v}^{\mathrm{sh}}, T)$,
	which is zero by \cite[Chapter III, Proposition 4.1]{Mil06}
	(where we used the assumption that $\kappa(v)$ is perfect).
	For
		$
			\operatorname{\mathbf{Ext}}^{n + 1}_{X_{\mathrm{sm}}}(
				i_{v, \ast} F_{v},
				\mathbf{G}_{m}
			)
		$,
	we may replace $X$ by the strict henselian localization at $v$
	and $F_{v}$ by $\mathbb{Z}$.
	Then
		\[
				\operatorname{\mathbf{Ext}}^{n + 1}_{X_{\mathrm{sm}}}(
					i_{v, \ast} \mathbb{Z},
					\mathbf{G}_{m}
				)
			\cong
				R^{n} j_{U, \ast} \mathbf{G}_{m},
		\]
	which is zero over $X_{\mathrm{et}}$ by the same reason.
\end{proof}

We need to pass on the ambiguity about mapping cones
``to the second variable'':

\begin{lem} \label{0018}
	Let $\mathcal{F} \colon \mathcal{A} \times \mathcal{A}' \to \mathcal{A}''$
	be a left exact bifunctor on abelian categories.
	Let $A' \in \mathcal{A}'$.
	Assume $\mathcal{A}'$ has enough injectives and
	$\mathcal{F}(\;\cdot\;, I') \colon \mathcal{A} \to \mathcal{A}''$
	is an exact functor
	for any injective object $I' \in \mathcal{A}'$
	(so that $\mathcal{F}$ admits a right derived functor
		$
				R \mathcal{F}
			\colon
				D^{+}(\mathcal{A}) \times D^{+}(\mathcal{A}')
			\to
				D^{+}(\mathcal{A}'')
		$
	by \cite[Corollary 13.4.5]{KS06}).
	\begin{enumerate}
		\item \label{0014}
			Let $0 \to A \to B \to C \to 0$ be an exact sequence in $\mathcal{A}$.
			Assume that $\mathcal{F}(A, A') = 0$.
			Then there exists a canonical morphism
				\begin{equation} \label{0013}
						[\mathcal{F}(C, A') \to R^{1} \mathcal{F}(A, A')][-1]
					\to
						\tau_{\le 1} R \mathcal{F}(B, A')
				\end{equation}
			in $D(\mathcal{A}'')$.
			It is the natural isomorphism
				\[
						\operatorname{Ker}(
							\mathcal{F}(C, A') \to R^{1} \mathcal{F}(A, A')
						)
					\cong
						\mathcal{F}(B, A')
				\]
			in $H^{0}$ and the natural injection
				\[
						\operatorname{Coker}(
							\mathcal{F}(C, A') \to R^{1} \mathcal{F}(A, A')
						)
					\hookrightarrow
						R^{1} \mathcal{F}(B, A')
				\]
			in $H^{1}$.
		\item \label{0015}
			Let $\{0 \to A_{\lambda} \to B \to C_{\lambda} \to 0\}_{\lambda \in \Lambda}$
			be a filtered direct system of exact sequences in $\mathcal{A}$
			with constant middle term.
			Assume that $\mathcal{A}''$ has exact filtered direct limits
			and $\mathcal{F}(A_{\lambda}, A') = 0$ for all $\lambda \in \Lambda$.
			Then there exists a canonical morphism
				\[
						\left[
								\varinjlim_{\lambda}
									\mathcal{F}(C_{\lambda}, A')
							\to
								\varinjlim_{\lambda}
									R^{1} \mathcal{F}(A_{\lambda}, A')
						\right][-1]
					\to
						\tau_{\le 1} R \mathcal{F}(B, A')
				\]
			in $D(\mathcal{A}'')$.
			This morphism and the morphism \eqref{0013} for each
			$0 \to A_{\mu} \to B \to C_{\mu} \to 0$
			fit in the commutative diagram
				\[
					\begin{CD}
							[\mathcal{F}(C_{\mu}, A') \to R^{1} \mathcal{F}(A_{\mu}, A')][-1]
						@>>>
							\tau_{\le 1} R \mathcal{F}(B, A')
						\\ @VVV @| \\
							\left[
									\varinjlim_{\lambda}
										\mathcal{F}(C_{\lambda}, A')
								\to
									\varinjlim_{\lambda}
										R^{1} \mathcal{F}(A_{\lambda}, A')
							\right][-1]
						@>>>
							\tau_{\le 1} R \mathcal{F}(B, A').
					\end{CD}
				\]
	\end{enumerate}
\end{lem}

\begin{proof}
	\eqref{0014}
	Let $\mathrm{Ch}^{+}(\mathcal{A})$ be
	the category of bounded below complexes in $\mathcal{A}$
	(\emph{not} up to homotopy!)
	and $\mathrm{Ch}^{+}(\mathcal{A}')$
	and $\mathrm{Ch}^{+}(\mathcal{A}'')$ similarly.
	First recall from \cite[Corollary 13.4.5]{KS06} that
	for any $D \in \mathrm{Ch}^{+}(\mathcal{A})$
	and $D' \in \mathrm{Ch}^{+}(\mathcal{A}')$,
	the object $R \mathcal{F}(D, D')$ can be calculated as
	the total complex of the double complex $\mathcal{F}(D, I')$,
	where $D' \to I'$ is a quasi-isomorphism
	to a bounded below injective complex in $\mathcal{A}'$.
	
	Let $B' \in \mathrm{Ch}^{+}(\mathcal{A}')$.
	Let $[A \to B] \in \mathrm{Ch}^{+}(\mathcal{A})$ be
	the mapping cone of $A \to B$.
	Let $h \colon [A \to B] \to B$ be the natural projection
	(which is a morphism of graded objects but not a morphism of complexes).
	Consider the mapping cone triangle
		\[
				B
			\to
				[A \to B]
			\to
				A[1]
			\to
				B[1]
		\]
	(in $\mathrm{Ch}^{+}(\mathcal{A})$).
	The composite $[A \to B] \to B[1]$ is null-homotopic
	and a null-homotopy is given by $h$.
	Consider the induced morphisms
		\[
				\tau_{\le 0} \mathcal{F}([A \to B], B')
			\to
				\tau_{\le 0} \mathcal{F}(A[1], B')
			\to
				\tau_{\le 0} \mathcal{F}(B[1], B')
		\]
	in $\mathrm{Ch}^{+}(\mathcal{A}'')$
	(where the functor $\tau_{\le 0}$ replaces
	the zero-th term by the kernel of the differential there
	and nullifies the positive degree terms).
	The composite
		\[
				\tau_{\le 0} \mathcal{F}([A \to B], B')
			\to
				\tau_{\le 0} \mathcal{F}(B[1], B')
		\]
	is null-homotopic and a null-homotopy is given by $\mathcal{F}(h, \mathrm{id})$.
	The choice of the null-homotopy $\mathcal{F}(h, \mathrm{id})$ defines a morphism
		\[
				\bigl[
						\tau_{\le 0} \mathcal{F}([A \to B], B')
					\to
						\tau_{\le 0} \mathcal{F}(A[1], B')
				\bigr]
			\to
				\tau_{\le 0} \mathcal{F}(B[1], B')
		\]
	in $\mathrm{Ch}^{+}(\mathcal{A}'')$
	by the proof of \cite[\href{https://stacks.math.columbia.edu/tag/041F}{Tag 014F}]{Sta22}
	or by \cite[Chapter I, Proposition 3.1.3]{Ver96}.
	This is functorial in $B'$.
	Hence we have a natural transformation between two functors
		$
				\mathrm{Ch}^{+}(\mathcal{A}')
			\rightrightarrows
				\mathrm{Ch}^{+}(\mathcal{A}'')
		$
	of additive categories with translation \cite[Definition 10.1.1]{KS06}
		\begin{equation} \label{0016}
				\bigl[
						\tau_{\le 0} \mathcal{F}([A \to B], \;\cdot\;)
					\to
						\tau_{\le 0} \mathcal{F}(A[1], \;\cdot\;)
				\bigr]
			\to
				\tau_{\le 0} \mathcal{F}(B[1], \;\cdot\;).
		\end{equation}
	These two functors induce functors
	$K^{+}(\mathcal{A}') \rightrightarrows K^{+}(\mathcal{A}'')$
	of additive categories with translation
	(where $K^{+}$ denotes the homotopy category)
	and the above gives a natural transformation between these two functors.
	Let $A' \to I'$ be a quasi-isomorphism
	to an injective complex concentrated in non-negative degrees
	(which is unique up to unique isomorphism in $K^{+}(\mathcal{A}')$).
	Then the above gives a morphism
		\[
				\bigl[
						\tau_{\le 0} \mathcal{F}([A \to B], I')
					\to
						\tau_{\le 0} \mathcal{F}(A[1], I')
				\bigr]
			\to
				\tau_{\le 0} \mathcal{F}(B[1], I')
		\]
	in $K^{+}(\mathcal{A}'')$ and hence in $D^{+}(\mathcal{A}'')$.
	We have canonical isomorphisms
		\begin{gather*}
					\tau_{\le 0} \mathcal{F}([A \to B], I')
				\cong
					\mathcal{F}(C, A'),
			\\
					\tau_{\le 0} \mathcal{F}(A[1], I')
				\cong
					R^{1} \mathcal{F}(A, A'),
			\\
					\tau_{\le 0} \mathcal{F}(B[1], I')
				\cong
					\tau_{\le 0} R \mathcal{F}(B[1], A')
				\cong
					(\tau_{\le 1} R \mathcal{F}(B, A'))[1]
		\end{gather*}
	in $D^{+}(\mathcal{A}'')$
	(where the second isomorphism used the assumption $\mathcal{F}(A, A') = 0$).
	Substituting these, we get the desired morphism.
	The required properties on $H^{0}$ and $H^{1}$ follow from the construction.
	
	\eqref{0015}
	The morphism \eqref{0016} and its functoriality give a commutative diagram
		\[
			\begin{CD}
					\bigl[
							\tau_{\le 0} \mathcal{F}([A_{\mu} \to B], \;\cdot\;)
						\to
							\tau_{\le 0} \mathcal{F}(A_{\mu}[1], \;\cdot\;)
					\bigr]
				@>>>
					\tau_{\le 0} \mathcal{F}(B[1], \;\cdot\;)
				\\ @VVV @| \\
					\left[
							\tau_{\le 0}
							\varinjlim_{\lambda}
							\mathcal{F}([A_{\lambda} \to B], \;\cdot\;)
						\to
							\tau_{\le 0}
							\varinjlim_{\lambda}
							\mathcal{F}(A_{\lambda}[1], \;\cdot\;)
					\right]
				@>>>
					\tau_{\le 0} \mathcal{F}(B[1], \;\cdot\;)
			\end{CD}
		\]
	between functors 
	$\mathrm{Ch}^{+}(\mathcal{A}') \to \mathrm{Ch}^{+}(\mathcal{A}'')$
	of additive categories with translation.
	The rest of the arguments goes similarly
	to the proof of the previous statement.
\end{proof}


\subsection{Duals of torsion-free constructible tori}

Now we will go in the opposite direction:
from torsion-free constructible tori
to torsion-free $\mathbb{Z}$-constructible sheaves.

\begin{prop} \label{0025}
	Let
		$
				\mathscr{T}
			=
				[
						\mathcal{T}
					\to
						\bigoplus_{v \in X_{0}}
							i_{v, \ast} E_{v}
				][-1]
		$
	be a torsion-free constructible torus over $X$.
	Then
		$
			\operatorname{\mathbf{Ext}}^{0}_{X_{\mathrm{sm}}}(
				\mathscr{T},
				\mathbf{G}_{m}
			)
		$
	is a torsion-free $\mathbb{Z}$-constructible sheaf over $X$
	and
		$
				\operatorname{\mathbf{Ext}}^{n}_{X_{\mathrm{sm}}}(
					\mathscr{T},
					\mathbf{G}_{m}
				)
			=
				0
		$
	for $n < 0$.
\end{prop}

\begin{proof}
	Assume that $\mathscr{T}$ has good reduction
	over a dense open subscheme $U \subset X$.
	Let
		$
				\mathscr{T}^{(U)}
			=
				[
						\mathcal{T}^{(U)}
					\to
						\bigoplus_{v \in X \setminus U}
							i_{v, \ast} E_{v}
				][-1]
		$
	be the $U$-model of $\mathscr{T}$.
	Let $F_{U}$ be the Cartier dual of $\mathcal{T}_{U}^{0}$.
	For each $v \in X_{0}$,
	let $F_{v}$ be the $\mathbb{Z}$-dual of $E_{v}$.
	Since the statement is invariant under quasi-isomorphism
	in $\operatorname{Ab}(X_{\mathrm{sm}})$
	and $\mathscr{T}^{(U)}$ is quasi-isomorphic to $\mathscr{T}$,
	we may replace $\mathscr{T}$ by $\mathscr{T}^{(U)}$.
	We have a distinguished triangle
	\adjustbox{center}{\resizebox{\textwidth}{!}{
		$
				R \operatorname{\mathbf{Hom}}_{X_{\mathrm{sm}}} \biggl(
					\displaystyle\bigoplus_{v \in X \setminus U}
						i_{v, \ast} E_{v},
					\mathbf{G}_{m}
				\biggr)
			\to
				R \operatorname{\mathbf{Hom}}_{X_{\mathrm{sm}}} \left(
					\mathcal{T}^{(U)},
					\mathbf{G}_{m}
				\right)
			\to
				R \operatorname{\mathbf{Hom}}_{X_{\mathrm{sm}}} \left(
					\mathscr{T},
					\mathbf{G}_{m}
				\right)
		$
	}}
	in $D(X_{\mathrm{sm}})$
	(where the connecting morphism of the triangle
	is not denoted here but assumed).
	By \cite[Theorem B (5)]{Suz19}, we have a canonical isomorphism
		\begin{equation} \label{0029}
				j_{U, !} F_{U}
			\overset{\sim}{\to}
				\tau_{\le 1}
				R \operatorname{\mathbf{Hom}}_{X_{\mathrm{sm}}} \left(
					\mathcal{T}^{(U)},
					\mathbf{G}_{m}
				\right)
		\end{equation}
	(again we used the assumption that
	the residue fields of $X$ at closed points are perfect).
	In particular, these isomorphic objects have trivial $H^{1}$.
	With Proposition \ref{0006} \eqref{0005},
	we know that
		$
				\operatorname{\mathbf{Ext}}^{n}_{X_{\mathrm{sm}}}(
					\mathscr{T}, \mathbf{G}_{m}
				)
			=
				0
		$
	for $n < 0$
	and we have an exact sequence
		\[
				0
			\to
				j_{U, !} F_{U}
			\to
				\operatorname{\mathbf{Ext}}^{0}_{X_{\mathrm{sm}}}(
					\mathscr{T},
					\mathbf{G}_{m}
				)
			\to
				\bigoplus_{v \in X \setminus U}
					i_{v, \ast} F_{v}
			\to
				0,
		\]
	showing that the middle term is torsion-free $\mathbb{Z}$-constructible.
\end{proof}

\begin{defi}
	Let $\mathscr{T}$ be a torsion-free constructible torus over $X$.
	Define
		\[
				\mathscr{T}^{D}
			=
				\operatorname{\mathbf{Ext}}^{0}_{X_{\mathrm{sm}}}(
					\mathscr{T},
					\mathbf{G}_{m}
				).
		\]
	For a morphism $\mathscr{T} \to \mathscr{S}$
	of torsion-free constructible tori over $X$,
	define a morphism $\mathscr{S}^{D} \to \mathscr{T}^{D}$
	to be the morphism induced by the functoriality of
	$\operatorname{\mathbf{Ext}}^{0}_{X_{\mathrm{sm}}}$.
\end{defi}

Thus we have a contravariant functor
$\mathscr{T} \mapsto \mathscr{T}^{D}$
from the category of torsion-free constructible tori over $X$
to the category of torsion-free $\mathbb{Z}$-constructible sheaves over $X$.
It is an additive functor.
Here is a presentation of $\mathscr{T}^{D}$:

\begin{prop} \label{0027}
	Let
		$
				\mathscr{T}
			=
				[
						\mathcal{T}
					\to
						\bigoplus_{v \in X_{0}}
							i_{v, \ast} E_{v}
				][-1]
		$
	be a torsion-free constructible torus over $X$.
	Let $U \subset X$ be a dense open subscheme
	where $\mathscr{T}$ has good reduction.
	Let $F_{U}$ be the Cartier dual of $\mathcal{T}_{U}^{0}$.
	For each $v \in X_{0}$,
	let $F_{v}$ be the $\mathbb{Z}$-dual of $E_{v}$.
	Then we have a canonical exact sequence
		\[
				0
			\to
				j_{U, !} F_{U}
			\to
				\mathscr{T}^{D}
			\to
				\bigoplus_{v \in X \setminus U}
					i_{v, \ast} F_{v}
			\to
				0.
		\]
	In particular, $\mathscr{T}^{D}$ is locally constant over $U$.
\end{prop}

\begin{proof}
	This was proved during the proof of Proposition \ref{0025}.
\end{proof}

We give an example of a dual:

\begin{prop} \label{0046}
	Let $\mathcal{T}$ and $\mathscr{T}$ be as in Example \ref{0045}.
	Let $F_{K}$ be the Cartier dual of $\mathcal{T}_{K}$.
	Then $\mathscr{T}^{D} \cong g_{\ast} F_{K}$.
\end{prop}

\begin{proof}
	Since $\mathscr{T}^{D}$ depends
	only on the image of $\mathscr{T}$ in $D(X_{\mathrm{sm}})$,
	we may identify $\mathscr{T}$ with $\mathcal{T}^{\tau}$.
	Consider the exact sequence
		\[
				0
			\to
				\mathcal{T}^{0}
			\to
				\mathcal{T}^{\tau}
			\to
				\bigoplus_{v \in X_{0}}
					i_{v, \ast}
					\pi_{0}(\mathcal{T}_{v})_{\mathrm{tor}}
			\to
				0
		\]
	in $\operatorname{Ab}(X_{\mathrm{sm}})$,
	where $\pi_{0}(\mathcal{T}_{v})_{\mathrm{tor}}$ is
	the torsion part of $\pi_{0}(\mathcal{T}_{v})$.
	Applying Proposition \ref{0006} \eqref{0005} to
	a torsion-free resolution of $\pi_{0}(\mathcal{T}_{v})_{\mathrm{tor}}$,
	we know that
		\[
				\operatorname{\mathbf{Ext}}^{n}_{X_{\mathrm{sm}}}(
					i_{v, \ast} \pi_{0}(\mathcal{T}_{v})_{\mathrm{tor}},
					\mathbf{G}_{m}
				)
			=
				0
		\]
	for $n \le 1$.
	Also, by \cite[Theorem B (5)]{Suz19},
	we have a canonical isomorphism
		\[
				\tau_{\le 1}
				R g_{\ast} F_{K}
			\overset{\sim}{\to}
				\tau_{\le 1}
				R \operatorname{\mathbf{Hom}}_{X_{\mathrm{sm}}}(
					\mathcal{T}^{0},
					\mathbf{G}_{m}
				).
		\]
	Hence the long exact sequence for
	$\operatorname{\mathbf{Ext}}^{\;\cdot\;}(\;\cdot\;, \mathbf{G}_{m})$
	gives
		$
				\operatorname{\mathbf{Hom}}_{X_{\mathrm{sm}}}(
					\mathcal{T}^{\tau},
					\mathbf{G}_{m}
				)
			\cong
				g_{\ast} F_{K}
		$.
\end{proof}


\subsection{Biduality}
\label{0037}

We will show that the category of torsion-free
$\mathbb{Z}$-constructible sheaves over $X$
and the category of torsion-free constructible tori over $X$
are contravariantly equivalent via the duality functors.

\begin{prop} \label{0034}
	Let $F$ be a torsion-free $\mathbb{Z}$-constructible sheaf over $X$.
	Then there exists a canonical functorial isomorphism
	$F \overset{\sim}{\to} F^{DD}$ ($:= (F^{D})^{D}$).
\end{prop}

\begin{proof}
	Consider the natural morphism
		\[
				F
			\to
				R \operatorname{\mathbf{Hom}}_{X_{\mathrm{sm}}} \bigl(
					R \operatorname{\mathbf{Hom}}_{X_{\mathrm{sm}}}(
						F,
						\mathbf{G}_{m}
					),
					\mathbf{G}_{m}
				\bigr)
		\]
	in $D(X_{\mathrm{sm}})$.
	Composing the morphism \eqref{0017} with it, we obtain a morphism
		\[
				F
			\to
				R \operatorname{\mathbf{Hom}}_{X_{\mathrm{sm}}}(
					F^{D},
					\mathbf{G}_{m}
				).
		\]
	Taking $H^{0}$, we obtain a morphism $F \to F^{DD}$.
	It is functorial.
	To show it is an isomorphism,
	let $U \subset X$ be a dense open subscheme where $F$ is locally constant.
	Denote the complex \eqref{0008} for $F$ by $F^{D_{U}}$
	(which is quasi-isomorphic to $F^{D}$
	in $\operatorname{Ab}(X_{\mathrm{sm}})$).
	Let $\mathcal{T}^{(U)}$ be the N\'eron model of the Cartier dual of $F_{U}$
	and $E_{v}$ the $\mathbb{Z}$-dual of $E_{v}$ for each $v \in X \setminus U$.
	We have a commutative diagram
		\begin{equation} \label{0030}
			\begin{CD}
					j_{U, !} F_{U}
				@>>>
					F
				@>>>
					\bigoplus_{v \in X \setminus U}
						i_{v, \ast} F_{v}
				\\ @VVV @VVV @VVV \\
					(j_{U, !} F_{U})^{DD}
				@>>>
					F^{DD}
				@>>>
					\bigoplus_{v \in X \setminus U}
						(i_{v, \ast} F_{v})^{DD}
			\end{CD}
		\end{equation}
	of $\mathbb{Z}$-constructible sheaves.
	We have \[(j_{U, !} F_{U})^{D_{U}} = \mathcal{T}^{(U)},
	\quad (i_{v, \ast} F_{v})^{D_{U}} = i_{v, \ast} E_{v}[-1]\]
	by definition.
	The sequence
		\[
				\bigoplus_{v \in X \setminus U}
					(i_{v, \ast} F_{v})^{D_{U}}
			\to
				F^{D_{U}}
			\to
				(j_{U, !} F_{U})^{D_{U}}
		\]
	can be identified with the term-wise exact sequence
		\begin{equation} \label{0028}
				0
			\to
				\bigoplus_{v \in X \setminus U}
					i_{v, \ast} E_{v}[-1]
			\to
				F^{D_{U}}
			=
				\left[
						\mathcal{T}^{(U)}
					\to
						\bigoplus_{v \in X \setminus U}
							i_{v, \ast} E_{v}
				\right][-1]
			\to
				\mathcal{T}^{(U)}
			\to
				0
		\end{equation}
	by construction.
	We have $(\mathcal{T}^{(U)})^{D} \cong j_{U, !} F_{U}$
	and $(i_{v, \ast} E_{v})^{D} \cong i_{v, \ast} F_{v}[-1]$
	by Proposition \ref{0027},
	where these upper $D$ mean
	$\operatorname{\mathbf{Ext}}^{0}_{X_{\mathrm{sm}}}$.
	Therefore by Proposition \ref{0025},
	the term-wise exact sequence \eqref{0028} induces an exact sequence
		\[
				0
			\to
				j_{U, !} F_{U}
			\to
				F^{DD}
			\to
				\bigoplus_{v \in X \setminus U}
					i_{v, \ast} F_{v}
		\]
	(which is actually right exact by \eqref{0029}).
	The morphisms
		\begin{gather*}
					j_{U, !} F_{U}
				\to
					(j_{U, !} F_{U})^{DD} \cong j_{U, !} F_{U},
			\\
					i_{v, \ast} F_{v}
				\to
					(i_{v, \ast} F_{v})^{DD} \cong i_{v, \ast} F_{v}
		\end{gather*}
	can be identified with the identity morphisms.
	Hence the diagram \eqref{0030} can be written as a commutative diagram
		\[
			\begin{CD}
					0
				@>>>
					j_{U, !} F_{U}
				@>>>
					F
				@>>>
					\bigoplus_{v \in X \setminus U}
						i_{v, \ast} F_{v}
				@>>>
					0
				\\ @. @| @VVV @| @. \\
					0
				@>>>
					j_{U, !} F_{U}
				@>>>
					F^{DD}
				@>>>
					\bigoplus_{v \in X \setminus U}
						i_{v, \ast} F_{v}
			\end{CD}
		\]
	with exact rows.
	A diagram chase then shows that
	$F \overset{\sim}{\to} F^{DD}$.
\end{proof}

\begin{prop} \label{0035}
	Let
		$
				\mathscr{T}
			=
				[
						\mathcal{T}
					\to
						\bigoplus_{v \in X_{0}}
							i_{v, \ast} E_{v}
				][-1]
		$
	be a torsion-free constructible torus over $X$.
	Then there exists a canonical functorial isomorphism
	$\mathscr{T} \overset{\sim}{\to} \mathscr{T}^{DD}$
	($:= (\mathscr{T}^{D})^{D}$).
\end{prop}

\begin{proof}
	Denote the given morphism
	$\mathcal{T} \to \bigoplus_{v \in X_{0}} i_{v, \ast} E_{v}$
	by $\varphi$.
	Let $U \subset X$ be a dense open subscheme
	where $\mathscr{T}$ has good reduction.
	Let
		$
				\mathscr{T}^{(U)}
			=
				[
						\mathcal{T}^{(U)}
					\to
						\bigoplus_{v \in X \setminus U}
							i_{v, \ast} E_{v}
				][-1]
		$
	be the $U$-model of $\mathscr{T}$.
	For a $\mathbb{Z}$-constructible sheaf $F$,
	let $F^{D_{U}}$ be the complex \eqref{0008}.
	It is enough to show that
	there exists a canonical functorial isomorphism
	$\mathscr{T}^{(U)} \overset{\sim}{\to} ((\mathscr{T}^{(U)})^{D})^{D_{U}}$
	of complexes of group algebraic spaces
	such that the diagram
		\begin{equation} \label{0033}
			\begin{CD}
					\mathscr{T}^{(U)}
				@> \sim >>
					((\mathscr{T}^{(U)})^{D})^{D_{U}}
				\\ @VVV @VVV \\
					\mathscr{T}^{(V)}
				@>> \sim >
					((\mathscr{T}^{(V)})^{D})^{D_{V}}
			\end{CD}
		\end{equation}
	commutes for any dense open $V \subset U$,
	where the right vertical morphism is the composite
		\[
				((\mathscr{T}^{(U)})^{D})^{D_{U}}
			\to
				((\mathscr{T}^{(U)})^{D})^{D_{V}}
			\to
				((\mathscr{T}^{(V)})^{D})^{D_{V}}.
		\]
	
	Denote
		$
				\operatorname{\mathbf{Hom}}
			=
				\operatorname{\mathbf{Hom}}_{X_{\mathrm{sm}}}
		$
	and
		$
				\operatorname{\mathbf{Ext}}^{n}
			=
				\operatorname{\mathbf{Ext}}^{n}_{X_{\mathrm{sm}}}
		$.
	Set $E^{(U)} = \bigoplus_{v \in X \setminus U} i_{v, \ast} E_{v}$
	and $F^{(U)} = \bigoplus_{v \in X \setminus U} i_{v, \ast} F_{v}$,
	where $F_{v}$ is the $\mathbb{Z}$-dual of $E_{v}$.
	Let $F = \mathscr{T}^{D} \cong (\mathscr{T}^{(U)})^{D}$.
	We have a commutative diagram
		\[
			\begin{CD}
					\mathcal{T}^{(U)}
				@>>>
					E^{(U)}
				\\ @VVV @VVV \\
					R \operatorname{\mathbf{Hom}}(
						R \operatorname{\mathbf{Hom}}(
							\mathcal{T}^{(U)},
							\mathbf{G}_{m}
						),
						\mathbf{G}_{m}
					)
				@>>>
					R \operatorname{\mathbf{Hom}}(
						R \operatorname{\mathbf{Hom}}(
							E^{(U)},
							\mathbf{G}_{m}
						),
						\mathbf{G}_{m}
					)
				\\ @VVV @VVV \\
					R \operatorname{\mathbf{Hom}}(
						\tau_{\le 1}
						R \operatorname{\mathbf{Hom}}(
							\mathcal{T}^{(U)},
							\mathbf{G}_{m}
						),
						\mathbf{G}_{m}
					)
				@>>>
					R \operatorname{\mathbf{Hom}}(
						\tau_{\le 1}
						R \operatorname{\mathbf{Hom}}(
							E^{(U)},
							\mathbf{G}_{m}
						),
						\mathbf{G}_{m}
					)
			\end{CD}
		\]
	in $D(X_{\mathrm{sm}})$.
	Since the upper row consists of objects concentrated in degree $0$,
	this diagram factors as
		\begin{equation} \label{0031}
			\adjustbox{center}{\resizebox{\textwidth}{!}{$\begin{CD}
					\mathcal{T}^{(U)}
				@>>>
					E^{(U)}
				\\ @VVV @VVV \\
					\tau_{\le 0}
					R \operatorname{\mathbf{Hom}}(
						\tau_{\le 1}
						R \operatorname{\mathbf{Hom}}(
							\mathcal{T}^{(U)},
							\mathbf{G}_{m}
						),
						\mathbf{G}_{m}
					)
				@>>>
					\tau_{\le 0}
					R \operatorname{\mathbf{Hom}}(
						\tau_{\le 1}
						R \operatorname{\mathbf{Hom}}(
							E^{(U)},
							\mathbf{G}_{m}
						),
						\mathbf{G}_{m}
					).
			\end{CD}$}}
		\end{equation}
	We have
		\begin{gather*}
					\tau_{\le 1}
					R \operatorname{\mathbf{Hom}}(
						\mathcal{T}^{(U)},
						\mathbf{G}_{m}
					)
				\cong
					j_{U, !} F_{U},
			\\
					\tau_{\le 1}
					R \operatorname{\mathbf{Hom}}(
						E^{(U)},
						\mathbf{G}_{m}
					)
				\cong
					F^{(U)}[-1]
		\end{gather*}
	by the proof of Proposition \ref{0025}.
	Therefore the lower row of the diagram \eqref{0031} can be written as
		\[
				\operatorname{\mathbf{Hom}}(
					j_{U, !} F_{U},
					\mathbf{G}_{m}
				)
			\to
				\operatorname{\mathbf{Ext}}^{1}(
					F^{(U)},
					\mathbf{G}_{m}
				)
		\]
	by Proposition \ref{0006} \eqref{0005}
	(all terms concentrated in degree zero).
	This morphism can be identified with the morphism \eqref{0032} for $F$
	since the connecting morphism
	$F^{(U)}[-1] \to j_{U, !} F_{U}$ comes from the morphism
		\[
				\varphi
			\colon
				\tau_{\le 1}
				R \operatorname{\mathbf{Hom}}(
					E^{(U)},
					\mathbf{G}_{m}
				)
			\to
				\tau_{\le 1}
				R \operatorname{\mathbf{Hom}}(
					\mathcal{T}^{(U)},
					\mathbf{G}_{m}
				).
		\]
	Thus \eqref{0031} gives a morphism
	$\mathscr{T}^{(U)} \to F^{D_{U}}$ of complexes.
	It is a term-wise isomorphism.
	This provides an isomorphism
	$\mathscr{T}^{(U)} \overset{\sim}{\to} ((\mathscr{T}^{(U)})^{D})^{D_{U}}$.
	It is functorial.
	The commutativity of the diagram \eqref{0033} is a diagram chase.
\end{proof}

\begin{thm} \label{0041}
	The two functors $D$ give
	mutually inverse contravariant equivalences
	between the category of torsion-free $\mathbb{Z}$-constructible sheaves over $X$
	and the category of torsion-free constructible tori over $X$.
\end{thm}

\begin{proof}
	This follows from Propositions \ref{0034} and \ref{0035}.
\end{proof}

\begin{ex} \label{0048}
	By Proposition \ref{0046},
	we have $(g_{\ast} F_{K})^{D} \cong \mathscr{T}$
	for the object $\mathscr{T}$ of Example \ref{0045}
	and the Cartier dual $F_{K}$ of $\mathcal{T}_{K}$.
\end{ex}


\subsection{Constructible tori (with torsion)}
\label{0038}

$\mathbb{Z}$-constructible sheaves may have torsion.
We will define a corresponding generalization
of torsion-free constructible tori.
Note that the opposite of the category of finitely generated abelian groups
is equivalent, via the functor
$R \operatorname{Hom}(\;\cdot\;, \mathbb{Z})$,
to the full subcategory of $D(\operatorname{Ab})$
consisting of objects
with finitely generated free $H^{0}$,
finite $H^{1}$ and $H^{n} = 0$ for $n \ne 0, 1$.
Hence we need to pass to a derived category
and then cut out some abelian category part.
As we already have the abelian category of $\mathbb{Z}$-constructible sheaves,
we use it as a model for our category of constructible tori
and do not have to go all the way through defining $t$-structures.
We will begin with defining a certain derived category
of the category of torsion-free constructible tori.

Let $\mathbb{Z} \mhyphen \mathrm{Con} / X$ be
the category of $\mathbb{Z}$-constructible sheaves over $X$
and $\mathbb{Z} \mhyphen \mathrm{Con}_{\mathrm{tf}} / X$
its full subcategory consisting of torsion-free sheaves.
Let $\mathrm{CT}_{\mathrm{tf}} / X$ be
the category of torsion-free constructible tori over $X$.
Since $\mathrm{CT}_{\mathrm{tf}} / X$ is an additive category,
we have its category of bounded complexes
$\mathrm{Ch}^{b}(\mathrm{CT}_{\mathrm{tf}} / X)$ and
its homotopy category
$K^{b}(\mathrm{CT}_{\mathrm{tf}} / X)$.

Let $\operatorname{Spec} K_{\mathrm{fppf}}$ be the fppf site of $K$.
Let $\operatorname{Ab}(K_{\mathrm{fppf}})$ be
the category of sheaves of abelian groups
on $\operatorname{Spec} K_{\mathrm{fppf}}$
and $D^{b}(K_{\mathrm{fppf}})$ its bounded derived category.

For a torsion-free $\mathbb{Z}$-constructible torus
	$
			\mathscr{T}
		=
			[
					\mathcal{T}
				\to
					\bigoplus_{v \in X_{0}}
						i_{v, \ast} E_{v}
			][-1]
	$
and a closed point $v \in X_{0}$, define
	\[
			\mathscr{T}_{K}
		=
			\mathcal{T}_{K},
		\quad
			\mathscr{T}_{v}
		=
			E_{v}.
	\]
They are tori over $K$ and a lattice over $v$, respectively.
View them as additive functors
	$
			\mathrm{CT}_{\mathrm{tf}} / X
		\to
			\operatorname{Ab}(K_{\mathrm{fppf}})
	$
and
	$
			\mathrm{CT}_{\mathrm{tf}} / X
		\to
			\operatorname{Ab}(v_{\mathrm{et}})
	$.
They extend to functors on the categories of (bounded) complexes.
Here is the definition of exact sequences and quasi-isomorphisms
in $\mathrm{CT}_{\mathrm{tf}} / X$:

\begin{defi} \mbox{}
	\begin{enumerate}
	\item
		An object $\mathscr{T}^{\bullet}$
		of $\mathrm{Ch}^{b}(\mathrm{CT}_{\mathrm{tf}} / X)$
		is said to be an \emph{exact complex}
		if $\mathscr{T}_{K}^{\bullet}$
		and $\mathscr{T}_{v}^{\bullet}$ for all $v \in X_{0}$
		are exact complexes
		in $\mathrm{Ch}^{b}(K_{\mathrm{fppf}})$
		and $\mathrm{Ch}^{b}(v_{\mathrm{et}})$, respectively.
		The sequence
		$\cdots \to \mathscr{T}^{n} \to \mathscr{T}^{n + 1} \to \cdots$
		is then said to be an \emph{exact sequence}
		in $\mathrm{CT}_{\mathrm{tf}} / X$.
	\item
		A morphism $\mathscr{T}^{\bullet} \to \mathscr{T}'^{\bullet}$
		in $\mathrm{Ch}^{b}(\mathrm{CT}_{\mathrm{tf}} / X)$
		is said to be a \emph{quasi-isomorphism}
		if its mapping cone is exact.
	\end{enumerate}
\end{defi}

Note that in this definition,
the N\'eron model part $\mathcal{T}^{\bullet}$ of $\mathscr{T}^{\bullet}$
is not required to be an exact complex
of smooth group schemes over $X$.
This definition matches with
that of torsion-free $\mathbb{Z}$-constructible sheaves:

\begin{prop} \label{0040}
	Let
		$
				F^{\bullet}
			\in
				\mathrm{Ch}^{b}(\mathbb{Z} \mhyphen \mathrm{Con}_{\mathrm{tf}} / X)
		$
	and
		$
				\mathscr{T}^{\bullet}
			\in
				\mathrm{Ch}^{b}(\mathrm{CT}_{\mathrm{tf}} / X)
		$
	be dual to each other
	(namely that one can be obtained from the other
	by applying the functor $D$ term-wise).
	Then $F^{\bullet}$ is an exact complex
	in the usual sense (in $\operatorname{Ab}(X_{\mathrm{et}})$)
	if and only if $\mathscr{T}^{\bullet}$ is an exact complex.
\end{prop}

\begin{proof}
	By definition,
	the objects $F_{K}^{\bullet}$ and $\mathscr{T}_{K}^{\bullet}$ are
	Cartier dual to each other
	and $F_{v}^{\bullet}$ and $\mathscr{T}_{v}^{\bullet}$
	for any $v \in X_{0}$ are $\mathbb{Z}$-duals to each other.
	Since $F^{\bullet}$ is exact if and only if
	$F_{K}^{\bullet}$ and $F_{v}^{\bullet}$ for all $v \in X_{0}$ are exact,
	the result follows.
\end{proof}

In particular, the exact complexes form
a full triangulated subcategory of the category
$K^{b}(\mathrm{CT}_{\mathrm{tf}} / X)$
closed under isomorphisms.
Therefore the following definition makes sense
by \cite[Definition 10.2.2 and Theorem 10.2.3]{KS06}:

\begin{defi}
	Define $D^{b}(\mathrm{CT}_{\mathrm{tf}} / X)$ to be
	the localization of $K^{b}(\mathrm{CT}_{\mathrm{tf}} / X)$
	at quasi-isomorphisms.
\end{defi}

\begin{prop}
	The term-wise application of the functors
	$\mathscr{T} \mapsto \mathscr{T}_{K}$ and
	$\mathscr{T} \mapsto \mathscr{T}_{v}$ for each $v \in X_{0}$ defines
	triangulated functors
		\[
				D^{b}(\mathrm{CT}_{\mathrm{tf}} / X)
			\to
				D^{b}(K_{\mathrm{fppf}})
			\quad \text{and} \quad
				D^{b}(\mathrm{CT}_{\mathrm{tf}} / X)
			\to
				D^{b}(v_{\mathrm{et}}),
		\]
	respectively.
\end{prop}

\begin{proof}
	These functors preserve quasi-isomorphisms by definition.
	Hence they factor through the localization at quasi-isomorphisms.
\end{proof}

We say that a morphism in
$K^{b}(\mathbb{Z} \mhyphen \mathrm{Con}_{\mathrm{tf}} / X)$
is a quasi-isomorphism
if it is so in
$K^{b}(\mathbb{Z} \mhyphen \mathrm{Con} / X)$
(or in $K^{b}(X_{\mathrm{et}})$).
Let
$D^{b}(\mathbb{Z} \mhyphen \mathrm{Con}_{\mathrm{tf}} / X)$
be the localization of
$K^{b}(\mathbb{Z} \mhyphen \mathrm{Con}_{\mathrm{tf}} / X)$
at quasi-isomorphisms.

\begin{prop} \label{0042}
	The duality functor $D$ and the natural inclusion functor give
	equivalences of triangulated categories
		\[
				D^{b}(\mathrm{CT}_{\mathrm{tf}} / X)
			\stackrel{D}{\simeq}
				D^{b}(
					\mathbb{Z} \mhyphen \mathrm{Con}_{\mathrm{tf}} / X
				)^{\mathrm{op}}
			\overset{\sim}{\to}
				D^{b}(\mathbb{Z} \mhyphen \mathrm{Con} / X)^{\mathrm{op}}.
		\]
\end{prop}

\begin{proof}
	The first equivalence follows from Propositions \ref{0041} and \ref{0040}.
	For the second, we know by \cite[\href{https://stacks.math.columbia.edu/tag/095N}{Tag 095N}]{Sta22} that
	any $\mathbb{Z}$-constructible sheaf over $X$ is
	a quotient of a torsion-free $\mathbb{Z}$-constructible sheaf.
	Also, a $\mathbb{Z}$-constructible subsheaf
	of a torsion-free $\mathbb{Z}$-constructible sheaf
	is torsion-free.
	Therefore we may apply \cite[Proposition 13.2.2 (ii)]{KS06},
	proving the desired equivalence.
\end{proof}

The proof shows that the inverse equivalence
	$
			D^{b}(\mathbb{Z} \mhyphen \mathrm{Con} / X)
		\overset{\sim}{\to}
			D^{b}(\mathbb{Z} \mhyphen \mathrm{Con}_{\mathrm{tf}} / X)
	$
sends a $\mathbb{Z}$-constructible sheaf $F$
to the complex $[F'' \to F']$,
where $F' \twoheadrightarrow F$ is any epimorphism
from a torsion-free $\mathbb{Z}$-constructible sheaf
and $F''$ is its kernel.
Now we define ``not necessarily torsion-free'' constructible tori:

\begin{defi}
	A \emph{constructible torus over $X$} is
	an object of $D^{b}(\mathrm{CT}_{\mathrm{tf}} / X)$
	isomorphic to a two-term complex
	$0 \to \mathscr{T}' \to \mathscr{T}'' \to 0$
	in $\mathrm{CT}_{\mathrm{tf}} / X$
	concentrated in degrees $0$ and $1$ such that
	$\mathscr{T}_{K}' \to \mathscr{T}_{K}''$ is faithfully flat
	and $\mathscr{T}_{v}' \to \mathscr{T}_{v}''$ has
	finite cokernel for all $v \in X_{0}$.
	The constructible tori over $X$ form
	a full subcategory of $D^{b}(\mathrm{CT}_{\mathrm{tf}} / X)$,
	which we denote by $\mathrm{CT} / X$.
\end{defi}

\begin{thm} \label{0043}
	The category of constructible tori over $X$ is
	contravariantly equivalent to
	the category of $\mathbb{Z}$-constructible sheaves over $X$
	under the equivalences in Proposition \ref{0042}.
	In particular, the category of constructible tori over $X$
	is an abelian category.
\end{thm}

\begin{proof}
	Let $F'' \to F'$ be a morphism of
	torsion-free $\mathbb{Z}$-constructible sheaves over $X$.
	Let $\mathscr{T}'$ and $\mathscr{T}''$ be
	the duals of $F'$ and $F''$, respectively.
	Then $F'' \to F'$ is a monomorphism if and only if
	$F''_{K} \to F'_{K}$ and $F''_{v} \to F'_{v}$ are
	monomorphisms for all $v \in X_{0}$
	if and only if
	$\mathscr{T}'_{K} \to \mathscr{T}''_{K}$ is faithfully flat
	and $\mathscr{T}'_{v} \to \mathscr{T}''_{v}$ has
	finite cokernel for all $v \in X_{0}$.
\end{proof}

Here is a characterization of constructible tori
among objects of $D^{b}(\mathrm{CT}_{\mathrm{tf}} / X)$:

\begin{prop}
	Let $\mathscr{T}^{\bullet}$ be
	an object of $D^{b}(\mathrm{CT}_{\mathrm{tf}} / X)$.
	Then $\mathscr{T}^{\bullet}$ is a constructible torus
	if and only if all of the following hold:
	\begin{enumerate}
		\item
			The object
			$\mathscr{T}_{K}^{\bullet} \in D^{b}(K_{\mathrm{fppf}})$
			is concentrated in degree zero.
		\item
			For any $v \in X_{0}$,
			the object
			$\mathscr{T}_{v}^{\bullet} \in D^{b}(v_{\mathrm{et}})$
			is concentrated in degrees $0$ and $1$
			with torsion-free $H^{0}$ and finite $H^{1}$.
	\end{enumerate}
\end{prop}

\begin{proof}
	Let
		$
				F^{\bullet}
			\in
				D^{b}(\mathbb{Z} \mhyphen \mathrm{Con}_{\mathrm{tf}} / X)
		$
	be the dual of $\mathscr{T}^{\bullet}$.
	Then the stated conditions can be translated to the following:
	\begin{enumerate}
		\item
			The object $F_{K}^{\bullet} \in D^{b}(K_{\mathrm{fppf}})$
			is concentrated in degree zero.
		\item
			For any $v \in X_{0}$,
			the object
			$F_{v}^{\bullet} \in D^{b}(v_{\mathrm{et}})$
			is concentrated in degree $0$.
	\end{enumerate}
	The conjunction of these conditions is equivalent to
	$F^{\bullet}$ being concentrated in degree zero
	whose cohomology is $\mathbb{Z}$-constructible.
	Thus the result follows from Proposition \ref{0043}.
\end{proof}

Since $\mathrm{CT} / X$ is an abelian category by Proposition \ref{0043},
its bounded derived category $D^{b}(\mathrm{CT} / X)$ makes sense.

\begin{prop} \label{0047}
	The natural functor
	$\mathrm{CT}_{\mathrm{tf}} / X \to \mathrm{CT} / X$
	is fully faithful.
	It induces an equivalence
		$
				D^{b}(\mathrm{CT}_{\mathrm{tf}} / X)
			\overset{\sim}{\to}
				D^{b}(\mathrm{CT} / X)
		$.
\end{prop}

\begin{proof}
	This follows from Propositions \ref{0042} and \ref{0043} by duality.
\end{proof}

\begin{prop}
	The natural functor
	$K^{b}(\mathrm{CT}_{\mathrm{tf}} / X) \to K^{b}(X_{\mathrm{et}})$
	induces a triangulated functor
	$D^{b}(\mathrm{CT} / X) \to D^{b}(X_{\mathrm{et}})$.
\end{prop}

\begin{proof}
	Consider the triangulated functors
		\[
				D^{b}(\mathrm{CT}_{\mathrm{tf}} / X)
			\stackrel{D}{\overset{\sim}{\to}}
				D^{b}(
					\mathbb{Z} \mhyphen \mathrm{Con}_{\mathrm{tf}} / X
				)^{\mathrm{op}}
			\hookrightarrow
				D^{b}(X_{\mathrm{et}})^{\mathrm{op}}
			\xrightarrow{
				R \operatorname{\mathbf{Hom}}_{X_{\mathrm{et}}}(
					\;\cdot\;,
					\mathbf{G}_{m}
				)
			}
				D(X_{\mathrm{et}}).
		\]
	Their composite is naturally identified with
	the one induced by the natural functor
	$\mathrm{CT}_{\mathrm{tf}} / X \to D^{b}(X_{\mathrm{et}})$
	by Proposition \ref{0026} \eqref{0021}.
	With Proposition \ref{0047}, the result follows.
\end{proof}

In particular, we can consider the \'etale cohomology complex
$R \Gamma(X_{\mathrm{et}}, \mathscr{T}^{\bullet})$
of any object
$\mathscr{T}^{\bullet} \in D^{b}(\mathrm{CT} / X)$.

\begin{rmk}
	A natural object of
	$D^{b}(\mathbb{Z} \mhyphen \mathrm{Con} / X)$
	not in $\mathbb{Z} \mhyphen \mathrm{Con} / X$
	is the N\'eron model $\mathcal{N}(F_{K})$ of a lattice $F_{K}$ over $K$
	in the sense of \cite[Theorem A]{Suz19},
	which is $\tau_{\le 1} R g_{\ast} F_{K}$
	as an object of $D(X_{\mathrm{sm}})$.
	Its dual in the sense of \cite[Theorem A (5)]{Suz19} is
	given by the connected N\'eron model
	$\mathcal{T}^{0}$ of the Cartier dual of $F_{K}$.
	On the other hand,
	our duality functor $D$ and Proposition \ref{0042} give
	an object $\mathcal{N}(F_{K})^{D}$ of $D^{b}(\mathrm{CT} / X)$
	not in $\mathrm{CT} / X$.
	Thus $\mathcal{N}(F_{K})^{D}$ can be thought of as $\mathcal{T}^{0}$
	``viewed as an object of $D^{b}(\mathrm{CT} / X)$.''
\end{rmk}

\subsection{Multiplicativity of the determinant of the Lie algebra for the Néron model of a torus}

In this section, we investigate the multiplicativity property of the determinant of the Lie algebra for Néron models of tori.

If $G$ is a group scheme over a scheme $S$, we denote $\scrLie(G)$ its Lie algebra sheaf \cite[Exp. II]{SGA3T1}; it is a sheaf of $\cal{O}_S$-modules on $S$. We also denote $\Lie(G)$ its global sections. Let $X$ be an irreducible Dedekind scheme (a noetherian regular scheme of dimension $\le 1$) with perfect residue fields, and denote by $K$ its function field. Then every torus $T$ over $K$ has a locally of finite type Néron model $\mathcal{N}(G)$ over $X$\cite[Chapter 10]{Bosch1990}. We begin with the following theorem, which is a rather immediate consequence of Chai--Yu's base change conductor formula for tori \cite{Chai2001}:
\begin{thm}[{\cite[Question 8.1]{Chai2000}}]\label{mult_det_lie_OX}
	Let
	\[
	0\to T' \to T \to T'' \to 0
	\]
	be a short exact sequence of $K$-tori. There is a multiplicative identity
	\[
	\det_{\cal{O}_X}\scrLie{} \cal{N}(T')\otimes_{\cal{O}_X} \det_{\cal{O}_X}\scrLie{}\cal{N}(T'') \xrightarrow{\simeq} \det_{\cal{O}_X}\scrLie{}\cal{N}(T),
	\]
	compatible under pullback to $K$ with the canonical isomorphism $\det_{K}\Lie(T')\otimes_K \det_{K}\Lie(T'') \xrightarrow{\simeq} \det_{K}\Lie(T)$.
\end{thm}
\begin{rmk}
	This is an identity of graded lines: since $X$ is connected, the compatibility of ranks can be checked after pulling back to $K$.
\end{rmk}
\begin{proof}
	It suffices to check this after localizing at a closed point so we reduce to the local case. Let $a\in K^\times$ be such that the image of $\det_{\cal{O}_K}\scrLie{}\cal{N}(T')\otimes_{\cal{O}_K} \det_{\cal{O}_K}\scrLie{}\cal{N}(T'')$ under the canonical isomorphism is $a\det_{\cal{O}_K}\scrLie{}\cal{N}(T)$; we want to show that $a\in \cal{O}_K^\times$. Let $L/K$ be a finite extension such that all the considered tori split over $L$. Let us check that $\det\scrLie{}\cal{N}(T'_L)\otimes \det\scrLie{}\cal{N}(T''_L) \xrightarrow{\simeq} \det\scrLie{}\cal{N}(T_L)$. For a split torus $S$ on $L$, we have $\scrLie{}\cal{N}(S)=\Hom_{\mathrm{Ab}}(X^\ast(S),\cal{O}_L)$ \cite[A1.7]{Chai2001} and the latter functor is exact in $S$ since the character group $X^\ast(T_L)$ is torsion-free, making the $\Ext^1$ vanish. We thus obtain the desired identity by taking determinants in the short exact sequence of Lie algebras of Néron models corresponding to the short exact sequence of split tori $0 \to T'_L \to T_L \to T''_L\to 0$. Recall Chai--Yu's base change conductor 
	\begin{align*}
		c(T)&:=\frac 1 {e_{L/K}} \len_{\cal{O}_L}\left(\frac{\scrLie{} \cal{N}(T_L)}{\scrLie{} \cal{N}(T)\otimes_{\cal{O}_K}\cal{O}_L}\right)\\
		&= \frac 1 {e_{L/K}} \len_{\cal{O}_L}\left(\frac{\det_{\cal{O}_L} \scrLie{} \cal{N}(T_L)}{(\det_{\cal{O}_K} \scrLie{} \cal{N}(T)) \otimes \cal{O}_L}\right).
	\end{align*}
	By Chai-Yu's formula expressing the base change conductor as one half of the Artin conductor of the cocharacter module \cite{Chai2001}, we have $c(T)=c(T')+c(T'')$ \cite[2.4.2]{Cluckers2013}. Consider the following diagram of sublattices
	\[\adjustbox{center}{\resizebox{\textwidth}{!}{\begin{tikzcd}[ampersand replacement=\&]
		{\det_{\cal{O}_K} \scrLie{} \cal{N}(T') \otimes_{\cal{O}_K} \det_{\cal{O}_K} \scrLie{} \cal{N}(T'') \otimes_{\cal{O}_K} \cal{O}_L}\& {\det_{\cal{O}_K} \scrLie{} \cal{N}(T') \otimes_{\cal{O}_K} \cal{O}_L} \\
		{\det_{\cal{O}_L} \scrLie{} \cal{N}(T'_L)\otimes_{\cal{O}_L} \det_{\cal{O}_L} \scrLie{} \cal{N}(T_L'')} \& {\det_{\cal{O}_L} \scrLie{} \cal{N}(T_L)} \\
		{\det_{L}\scrLie{}(T'_L)\otimes_L \det_{L}\scrLie{}(T''_L)} \& {\det_{L}\scrLie{}(T_L)}
		\arrow["\subset"{description}, sloped, draw=none, from=1-2, to=2-2]
		\arrow["\simeq", from=2-1, to=2-2]
		\arrow["\subset"{description}, sloped, draw=none, from=1-1, to=2-1]
		\arrow["\subset"{description}, sloped, draw=none, from=2-1, to=3-1]
		\arrow["\subset"{description}, sloped, draw=none, from=2-2, to=3-2]
		\arrow["\simeq", from=3-1, to=3-2]
	\end{tikzcd}}}\]
	Computing lengths, we find\footnote{Using some light abuse of notations, where we denote $\len(M/N):=-\len(N/M)$ if $M\subset N$ for two $\cal{O}_L$ lattices $M$ and $N$ inside $L$}
	\begin{align*}
		v_L(a)=&\len\left(\frac{a\det_{\cal{O}_K}\scrLie{}\cal{N}(T)\otimes_{\cal{O}_K}\cal{O}_L}{\det_{\cal{O}_K}\scrLie{}\cal{N}(T)\otimes_{\cal{O}_K}\cal{O}_L}\right)\\&=\len\left(\frac{a\det_{\cal{O}_K}\scrLie{}\cal{N}(T)\otimes_{\cal{O}_K}\cal{O}_L}{\det_{\cal{O}_L} \scrLie{} \cal{N}(T_L)}\right) + e_{L/K}c(T)\\ &= \len\left(\frac	{\det_{\cal{O}_K} \scrLie{} \cal{N}(T') \otimes_{\cal{O}_K} \det_{\cal{O}_K} \scrLie{} \cal{N}(T'') \otimes_{\cal{O}_K} \cal{O}_L}{\det_{\cal{O}_L} \scrLie{} \cal{N}(T'_L)\otimes_{\cal{O}_L} \det_{\cal{O}_L} \scrLie{} \cal{N}(T_L'')}\right) + e_{L/K}c(T)\\
		&= e_{L/K}(- c(T') - c(T'') + c(T))\\
		&=0.
	\end{align*}
\end{proof}

\begin{cor}\label{mult_det_lie_neron}
	Suppose that $X$ is proper over $\Spec(\ZZ)$, so that $K$ is a global field. Let
	\[
	0\to T' \to T \to T'' \to 0
	\]
	be a short exact sequence of tori over $K$. We have a multiplicative identity
	\[
	\det_{\ZZ}R\Gamma(X,\scrLie{}\cal{N}(T'))\otimes \det_{\ZZ}R\Gamma(X,\scrLie{}\cal{N}(T'')) \xrightarrow{\simeq} \det_{\ZZ} R\Gamma(X,\scrLie{}\cal{N}(T)),
	\]
	compatible with $\det_{\QQ}\Lie_K(T')\otimes_\QQ \det_{\QQ}\Lie_K(T'') \xrightarrow{\simeq} \det_{\QQ}\Lie_K(T)$ in the number field case, and with the canonical trivializations $\bigl(\det_{\ZZ} R\Gamma(X,\scrLie{}\cal{N}(T))\bigr)\otimes\QQ \xrightarrow{\simeq} \QQ$ in the function field case.
\end{cor}
\begin{rmk}
	If $K$ is a number field then we see $R\Gamma(X,-)$ as the forgetful functor $D^b(X)=D^b(\cal{O}_K)\to D^b(\ZZ)$.
\end{rmk}

\begin{proof}
	Suppose that $K$ is a number field. For a commutative ring $R$, denote $\cal{P}_R$ the Picard groupoid of graded $R$-lines. By $K$-theory computations in degree $\leq 1$, the determinant functors $\det_{R} :D^{perf}(R) \to \cal{P}_{R}$ for $R=\cal{O}_K,K,\ZZ,\QQ$ are the universal determinant functors \cite[§ 2.5]{Burns2003}. Composition of any determinant functor on $D^{perf}(\ZZ)$ with the forgetful functor $U=R\Gamma(X,-):D^{perf}(\cal{O}_K)\to D^{perf}(\ZZ)$ induces a determinant functor, so we get an induced map of Picard groupoids $\widetilde{U}:\cal{P}_{\cal{O}_K}\to \cal{P}_\ZZ$ such that $\widetilde{U}\circ \det_{\cal{O}_K}=\det_\ZZ \circ U$. We can thus apply the functor of Picard groupoids $\widetilde{U}$ to the multiplicative identity of the previous theorem to obtain the identity we want (implicit is the compatibility with base change and the analogous forgetful functor $D^b(K)\to D^b(\QQ)$).
	
	If $K$ is a function field, the corollary amounts to showing the equality of Euler characteristics of coherent sheaves
	\[
	\chi_X(\scrLie{}\cal{N}(T))=\chi_X(\scrLie{}\cal{N}(T'))+\chi_X(\scrLie{}\cal{N}(T''))
	\]
	By Riemann-Roch, we have
	\[
	\chi_X(\scrLie{}\cal{N}(T))=(1-g)\dim T + \deg(\det_X \scrLie{}\cal{N}(T))
	\]
	whence we conclude from the previous theorem.
\end{proof}
\begin{rmk}
	Note that the argument for the number field case does not go through for the function field case because the usual determinant is \emph{not} the universal determinant. Indeed we have $K_1(X)=(\cal{O}_X(X)^\times)^2$ (\cite[2.5]{Dwyer1997}) while $\pi_1(\cal{P}_X)=\cal{O}_X(X)^\times$.
\end{rmk}


\subsection{The functor \texorpdfstring{$\det\Lie$}{detLie} for constructible tori}
\label{0039}

We will define the Lie algebra of a torsion-free constructible torus.
For more general constructible tori and
even more general objects of $D^{b}(\mathrm{CT} / X)$,
unfortunately their Lie algebras do not make sense.
Only the determinants of Lie algebras do.

For a presheaf of groups $G$ on the category of $X$-schemes,
let $\operatorname{Lie}(G) = \operatorname{Lie}(G / X)$
be its Lie algebra presheaf
(\cite[Expos\'e II, Section 3.8]{DG70b},
\cite[Chapter II, Section 4, No.\ 1.2]{DG70a}).
If $G$ is a smooth group scheme
(or more generally, a smooth group algebraic space),
then $\operatorname{Lie}(G)$ is a locally free $\mathcal{O}_{X}$-module
dual to the pullback of $\Omega_{G / X}^{1}$
by the identity section $X \to G$
(\cite[Expos\'e II, Section 4.11]{DG70b},
\cite[Chapter II, Section 4, Proposition 4.8]{DG70a}).

\begin{defi}
	Let
		$
				\mathscr{T}
			=
				[
						\mathcal{T}
					\to
						\bigoplus_{v \in X_{0}}
							i_{v, \ast} E_{v}
				][-1]
		$
	be a torsion-free constructible torus over $X$.
	Define its Lie algebra by applying
	$\operatorname{Lie}$ term-wise:
		\[
				\operatorname{Lie}(\mathscr{T})
			=
				\left[
						\operatorname{Lie}(\mathcal{T})
					\to
						\bigoplus_{v \in X_{0}}
							\operatorname{Lie}(i_{v, \ast} E_{v})
				\right][-1]
			\cong
				\operatorname{Lie}(\mathcal{T}).
		\]
\end{defi}

Here we used the facts that
the Lie algebra functor is left exact
and annihilates formally \'etale presheaves by definition.
We thus have a functor $\operatorname{Lie}$
from the category of torsion-free constructible tori over $X$
to the category of locally free $\mathcal{O}_{X}$-modules of finite rank.
Taking the highest exterior power $\det_{\mathcal{O}_{X}}$,
we obtain an invertible $\mathcal{O}_{X}$-module
$\det_{\mathcal{O}_{X}} \operatorname{Lie}(\mathscr{T})$.

\begin{prop}\label{prop:extension_detLie}
    The assignment $\mathscr{T} \mapsto \det_{\mathcal{O}_{X}} \operatorname{Lie}(\mathscr{T})$, resp. $\mathscr{T} \mapsto \det_{\ZZ} \operatorname{Lie}(\mathscr{T})$ defines a determinant functor on the exact category of torsion-free constructible tori, valued in graded $\cal{O}_X$-lines (resp. graded $\ZZ$-lines). It extends uniquely to a determinant functor $\det\Lie$ (resp. $\Delta_X^{\mathrm{add}}$) on $D^b(\mathrm{CT}/X)$. In particular, we have a functor $\Delta_X^{\mathrm{add}}:D^b(\mathrm{CT}/X)_{\mathrm{iso}} \to \cal{P}_{\ZZ}$ from the maximal subgroupoid to the Picard groupoid of graded $\ZZ$-lines.
\end{prop}

\begin{proof}
We do the proof for $\det\Lie$, the proof is completely analogous for $\Delta^{\mathrm{add}}$. The axioms for determinants on exact categories are laid out in \cite[§1]{Knudsen2002}. Since $\det_{\cal{O}_X}$ is a determinant, to obtain a predeterminant structure the only thing to check is the multiplicativity with respect to short exact sequences. However, given a short exact sequence $0 \to \mathscr{T}' \to \mathscr{T} \to  \mathscr{T}'' \to 0$, the corresponding sequence $0\to T' \to T \to T'' \to 0$ on $K$ is exact; since $\Lie(\mathscr{T})= \Lie(\cal{T})$, the multiplicativity is \cref{mult_det_lie_OX}. Denote by $g:\Spec(K) \to X$ the inclusion of the generic point. The pullback $g^\ast$ of quasicoherent sheaves is exact and induces a pullback functor on graded lines, and we have a canonical base change isomorphism $g^\ast\circ \det_{\cal{O}_X}\simeq \det_K \circ g^\ast$. Moreover, $g^\ast$ is faithful; hence we can check that the compatibility, associativity and commutativity axioms hold after applying $g^\ast$. Since $g^\ast\det_{\mathcal{O}_{X}} \operatorname{Lie}(\mathscr{T})\simeq \det_K g^\ast \operatorname{Lie}(\mathscr{T})\simeq \det_K \Lie_K(\mathscr{T}_K)$, canonically, those axioms hold because $\mathscr{T} \mapsto \det_K \Lie_K(\mathscr{T}_K)$ is the composition of a determinant functor with an exact functor.

Since $K$-theory and determinant functors are insensitive to passing to the opposite category, we conclude with \cref{K-theory_torsionfree}.
\end{proof}

\begin{ex}
	On the other hand, the Lie algebra functor
	$\mathscr{T}^{\bullet} \mapsto \operatorname{Lie}(\mathscr{T}^{\bullet})$
	itself does not factor through $D^{b}(\mathrm{CT} / X)$.
	To see an example,
	let $X = \operatorname{Spec} \mathcal{O}_{K} = \operatorname{Spec} k[[t]]$
	with $k$ a perfect field of characteristic $2$ and with fraction field $K$.
	Let $\mathcal{O}_{L} = \mathcal{O}_{K}[x] / (x^2 + t x + t)$
	with fraction field $L$.
	Let $\operatorname{Res}_{L / K} \mathbf{G}_{m}$ be
	the Weil restriction of the multiplicative group.
	Let $T$ be the kernel of the norm map
	$\operatorname{Res}_{L / K} \mathbf{G}_{m} \to \mathbf{G}_{m}$; notice that it its character group is the sign representation on $\ZZ$.
	Let $\mathcal{T}$ be its N\'eron model over $\mathcal{O}_{K}$.
	Then we have an exact sequence
		\begin{equation} \label{0049}
				0
			\to
				\mathcal{T}
			\to
				\operatorname{Res}_{\mathcal{O}_{L} / \mathcal{O}_{K}}
					\mathbf{G}_{m}
			\to
				 \mathbf{G}_{m}
			\to
				0
		\end{equation}
	of torsion-free constructible tori over $\mathcal{O}_{K}$
	(where we are viewing relative torsion components of N\'eron models
	as torsion-free constructible tori using Example \ref{0045}).
	Hence this sequence viewed as an object of
	$D^{b}(\mathrm{CT}_{\mathrm{tf}} / X)$
	is zero.
	However, it is not an exact sequence in $\operatorname{Ab}(X_{\mathrm{sm}})$.
	In fact, if we apply Lie term-wise to \eqref{0049},
	we get the sequence
		$
				0
			\to
				\mathcal{O}_{K}
			\to
				\mathcal{O}_{K}^{2}
			\to
				\mathcal{O}_{K}
			\to
				0
		$,
	where the first map sends $a$ to $(t a, 0)$ and
	the second sends $(b, c)$ to $t c$.
	This is not an exact sequence
	(though its determinant is trivial as expected).
\end{ex}

%% file: specialValue.tex
\subsection{\texorpdfstring{$L$}{L}-functions for constructible tori}
\label{0044}

Let $X$ be a regular $1$-dimensional integral arithmetic scheme. For any scheme $S$, we will denote $\nu:\Sh(S_{\mathrm{proet}})\to \Sh(S_{\mathrm{et}})$ the natural morphism of topoi; the left adjoint $\nu^\ast$ is fully faithful \cite[lem. 5.1.2]{Bhatt2015}.  Denote by $-\widehat{\otimes}\QQ_l:=(R\lim (\nu^\ast(-)\otimes^L\ZZ/l^n\ZZ))\otimes \QQ$ the \emph{completed} tensor product of an étale sheaf with $\QQ_l$. If $\scr{T}$ is a constructible torus on $X$, its étale realization allows the definition of an $L$-function as in \cite{Morin2023b}:

\begin{defi}
	For each closed point $x$ of $X$, let $l_x$ be a prime number coprime to the residual characteristic at $x$ and $\varphi$ be the \emph{geometric} Frobenius in $G_x:=\mathrm{Gal}(\kappa(x)^{\mathrm{sep}}/\kappa(x))$. 
	Let $\scr{T}^\bullet \in D^b(\mathrm{CT}/X)$. We define the $L$-function of $\scr{T}^\bullet$ by the Euler product
	\[
	L_X(\scr{T},s)=\prod_{x\in X_0} \det\bigl (I-N(x)^{-s}\varphi|\bigl( i_x^\ast \scr{T}^\bullet_{\mathrm{et}}\bigr )\widehat{\otimes}\QQ_{l_x}\bigr )^{-1}
	\]
\end{defi}
We will see below a more explicit description, which implies in particular that this $L$-function doesn't depend on the choice of the prime numbers $(l_x)_{x\in X_0}$.

Let $T$ be torus over $K$ with character group $Y$. For any prime $l$, define the rational $l$-adic Tate module of $T$ as
\[
V_l(T):=(\lim T[l^n])\otimes \QQ
\]
It is a finite dimensional $l$-adic representation of $G_K$. We know that
\[
V_l(T) \simeq Y^\vee \otimes \QQ_l(1)
\]
as $l$-adic representations. 

\begin{defi}[{\cite[§ 8]{Fontaine1992}}]
	The $L$-function $L_K(T,s)$ of $T$ is the $L$-function of the $1$-motive $[0 \to T]$ over $K$, defined by the Euler product:
	\[
	L(T,s)=\prod_{v \text{ finite place of } K}\det(I-N(v)^{-s}\varphi|V_{l_v}(T)^{I_v})^{-1}
	\]
	where $l_v$ is a prime number coprime to the residual characteristic at $v$. 
\end{defi}
By the above, the $L$-function of $T$ is also the Artin $L$-function at $s+1$ of $Y\otimes \QQ$, and it doesn't depend on the choice of the family $(l_v)$.
\begin{rmk}
	This differs from the definition in \cite{Geisser2022}: it uses the $L$-function of the $1$-motive $[0\to T]$, involving $V_l(T)$, while Geisser--Suzuki use the Hasse-Weil $L$-function of $T$ involving $V_l(T)(-1)$.
\end{rmk}

 Let $i:x \to X$ be a closed point. Consider $K_x^h=\mathrm{Frac}(\cal{O}_x^h)$ the henselian local field at $x$. Fix an embedding $K_x^{h} \hookrightarrow K^{\mathrm{sep}}$. This determines an inertia group $I_x$ inside $G_K$ which is the absolute Galois group of $K_x^{sh}$. If $N$ is a discrete $G_x$-module with free of finite type underlying abelian group, or a rational or $l$-adic $G_x$-representation of finite dimension, denote
\[
Q_x(N,s):= \det(I-N(x)^{-s}\varphi| N)^{-1}
\]
with $\varphi$ the geometric Frobenius at $x$.

We recall the following:
\begin{thm}[{\cite[§ 6]{Morin2023b}}]
	Let $\scr{T}$ be a torsion-free constructible torus on $X$, and denote by $T$ its pullback to $\Spec(K)$ and $Y$ the character group of $T$.
	\begin{enumerate}
	\item The local factor at $x$ of the $L$-function of $\scr{T}$ is
	\[
	L_x(\scr{T},s)=\frac{Q_x((Y^\vee)^{I_x},s)}{Q_x((Y^\vee)^{I_x},s+1)Q_x(\scr{T}_x,s)}
	\]
	\item If $\scr{T}$ has good reduction around the closed point $x$, then $I_x$ acts trivially on $Y$, $Y=\scr{T}_x^\vee$ and we find that the local factor at $x$ of $\scr{T}$ equals the inverse of the local factor at $x$ of the $L$-function of $T$:\footnote{Notice that any rational representation is self-dual}
	\[
	L_x(\scr{T},s)=L_x(T,s)^{-1}
	\]
	\item Each local factor is well-defined, independently of the choice of a prime number $l$ coprime to the residual characteristic.
	\item The $L$-function of $\scr{T}$ differs from $L(T,s)^{-1}$ by a finite number of factors and is thus well-defined.
	\item $L_X(\scr{T},s)$ is meromorphic.
	\item Suppose that $X$ is proper. We have $L_X((g_\ast Y)^D,s)=L(T,s)^{-1}$.
	\item Suppose that $X$ is proper and denote $\cal{T}^0$ the connected (lft) Néron model of $T$ on $X$. 
	\[
	L_X(\cal{T}^0,s)=L_K(T,s)^{-1}
	\]
\end{enumerate}
\end{thm}

\begin{prop}
	Let $0\to \scr{T}^{\bullet} \to \scr{T} \to \scr{T}'' \to 0$ be a short exact sequence of constructible tori on $X$. Then
	\[
	L_X(\scr{T},s)=L_X(\scr{T}',s)L_X(\scr{T}'',s)
	\]
	Let $\pi:Y\to X$ be a finite dominant morphism with $Y$ regular integral, and let $\scr{T}$ be a constructible torus on $Y$. We have\footnote{See the next section for the definition of $\pi_\ast$}
	\[
	L_X(\pi_\ast \scr{T},s)=L_Y(\scr{T},s)
	\]
\end{prop}


\subsection{The fundamental line and special values at $s=0$}

Let $X$ be a regular $1$-dimensional integral proper arithmetic scheme. Denote by $K$ the function field of $X$. Then $X$ is either the spectrum of a ring of integers in the number field $K$ or a proper smooth curve over a finite field with function field $K$.

We prove a Weil-étale special value formula for $L_{X}(\mathscr{T}^{\bullet}, s)$ at $s = 0$ for any object $\mathscr{T}^{\bullet} \in D^{b}(\mathrm{CT}/X)$, thus extending the results of \cite{Morin2023b} and \cite{Geisser2022} to all constructible tori on $X$.

\begin{defi}
    Let $\mathscr{T}^\bullet \in D^b(\mathrm{CT}/X)$. The étale realization of $\mathscr{T}^\bullet$ is of the form $R\sheafhom_{X_\mathrm{et}}((\mathscr{T}^\bullet)^D,\GG_m)$, so the cohomology with Betti compact support at infinity $R\Gamma_{\mathrm{B}}(X,\mathscr{T}^\bullet)$ is defined as in \cite[Def. 2.1]{Morin2023b} through the fiber sequence
    \[
    R\Gamma_{\mathrm{B}}(X,\mathscr{T}^\bullet)\to R\Gamma(X,\mathscr{T}^\bullet) \to \prod_{v\mid \infty} R\Hom_{G_v}(Y^\bullet,2i\pi\ZZ[1])
    \]
    In the function field case, it is equal to ordinary étale cohomology. We say that a constructible torus $\mathscr{T}$ is:
    \begin{itemize}
        \item \emph{red} if $H^1_{\mathrm{B}}(X,\mathscr{T})$ is finite;
        \item \emph{blue} if $H^2_{\mathrm{B}}(X,\mathscr{T})$ is finite.\footnote{Beware the shift in indexing conventions between this article and \cite{Morin2023b}; namely, we have $F^D=\mathscr{T}_{\mathrm{et}}[1]$ where $F^D$ denotes the dual of a $\ZZ$-constructible sheaf as in loc. cit.}
    \end{itemize}
    A blue-to-red morphism is a morphism either between blue constructible tori, red constructible tori, and from a blue to a red one. A blue-to-red short exact sequence is a short exact sequence with blue-to-red morphisms.
\end{defi}

It $\cal{T}_U$ is a torus on a dense open $U$ of $X$, then (the torsion-free constructible torus representing) the Néron model of $\cal{T}_U$ is red, and a connected Néron model is red\footnote{By \cite[thm. B]{Suz19}, the dual of the connected Néron model of a $K$-torus $T$ is the complex $\tau_{\leq 1}Rg_\ast Y$ where $Y$ is the character group of $T$. Hence a connected Néron model is in general not a torsion-free constructible torus, but lives in the derived category.}. A constructible torus with no generic part (equivalently such that the dual has finite support) is blue. In particular, any constructible torus $\scr{T}$ fits in a short exact sequence
\[
0\to \scr{T}'\to \scr{T}\to \scr{T}''\to 0
\]
where $\scr{T}'$ is blue and $\scr{T}''$ is red\footnote{This can for instance be seen on the level of duals, where any $\ZZ$-constructible sheaf $F$ fits in a short exact sequence $0\to j_!F_{\mid U} \to F \to \oplus_{x\notin U} i_{x,\ast} i_x^\ast F \to 0$ where $j:U\to X$ is a dense open such that $F_{\mid U}$ is locally constant.} 
In order to define an Euler characteristic on $D^b(\mathrm{CT}/X)$ it thus suffices to restrict ourselves to red and blue constructible tori.

\begin{defiprop}
    Let $\mathscr{T}$ be a red or blue constructible torus on $X$. Then its Weil-étale cohomology with Betti compact support at infinity $R\Gamma_{\mathrm{W,B}}(X,\scr{T})$ is defined as in \cite[§ 3.1]{Morin2023b} through a distinguished triangle coming from Artin-Verdier duality
    \[
    R\Hom(R\Gamma(X,\scr{T}^D),\QQ[-3]) \to R\Gamma_{\mathrm{B}}(X,\mathscr{T})\to R\Gamma_{\mathrm{W},\mathrm{B}}(X,\mathscr{T})
    \]
    It is a perfect complex defined uniquely up to unique isomorphism, it is functorial in blue-to-red morphisms and finite dominant pushforwards and gives a long exact cohomology sequence for blue-to-red short exact sequences.
\end{defiprop}
\begin{proof}
    See \cite[§3.1]{Morin2023b} for the number field case; all proofs are valid in the function field case.
\end{proof}

\begin{prop}
    Let $\scr{T}$ be a red or blue constructible torus. We have a canonical rational splitting:
    \[
    R\Gamma_{\mathrm{W,B}}(X,\scr{T})_\QQ = R\Hom(R\Gamma(X,\scr{T}^D),\QQ[-2])\oplus R\Gamma_{\mathrm{B}}(X,\scr{T})_\QQ
    \]
\end{prop}

\begin{proof}
    See \cite[prop. 3.12]{Morin2023b}.
\end{proof}

\begin{defi}
    Let $\scr{T}$ be a red or blue constructible torus on $X$. The fundamental line is the graded $\ZZ$-line
    \[
    \Delta_X(\scr{T}):=\bigl(\det_\ZZ R\Gamma_{\mathrm{W},\mathrm{B}}(X,\scr{T})\bigr)^{-1}\otimes_\ZZ \Delta_X^{\mathrm{add}}(\scr{T})
    \]
\end{defi}
When the $\scr{T}$ is tamely ramified, i.e. when $\scr{T}_K$ splits after a tamely ramified extension (for instance if $\scr{T}=\GG_m$), this is the inverse of $\Delta_X(F^D)$ as defined in \cite{Morin2023b}, because in the latter the definition of $F^D$ involves a shift.

Let $Y$ be a regular integral scheme with a finite dominant morphism $\pi:Y\to X$. Consider the functor $\mathrm{CT}_{\mathrm{tf}}/Y\to\mathrm{CT}_{\mathrm{tf}}/X$, $\scr{T}\mapsto \pi_\ast \scr{T}$ which sends $\scr{T}=[\mathcal{T}\to\bigoplus_{w \in Y_{0}}i_{w, \ast} E_{v}][-1]$ to
\[
\pi_\ast \scr{T}:=[ \mathrm{Res}_{Y/X}\mathcal{T} \to \bigoplus_{v \in X_{0}} 	i_{v, \ast} \bigoplus_{\pi(w)=v} \pi_{w,\mathrm{et},\ast}E_{w}][-1]
\]
where $\pi_w$ is the induced morphism $w\to \pi(w)$. Since Weil restriction represents $\pi_{\mathrm{fppf},\ast}$, we have $\pi_{\ast}\scr{T}=\pi_{\mathrm{sm},\ast} \scr{T}$ as chain complexes of sheaves on $X_{\mathrm{sm}}$. As Weil restriction commutes with the formation of Néron models \cite[Prop. 7.6/6]{Bosch1990}, $\pi_\ast(\scr{T})$ is the dual of $\pi_{\mathrm{et},\ast} (\scr{T}^D)$; in particular $\pi_\ast$ is exact so it extends to a functor $\pi_\ast:D^b(\mathrm{CT}/Y)\to D^b(\mathrm{CT}/X)$ which is just the opposite of $\pi_{\mathrm{et},\ast}$. Because the right adjoint $R\pi_{\mathrm{et}}^!$ of $\pi_{\et,\ast}$ satisfies $R\pi_{\et}^! \GG_m=\GG_m$, this also implies that the étale realization of $\pi_\ast$ is $\pi_{\mathrm{et},\ast}$.
\begin{prop}\label{compatibilty_pushforward}
Let $Y$ be a regular integral scheme with a finite dominant morphism $\pi:Y\to X$ and let $\scr{T}$ be a torsion-free constructible torus on $Y$. Then
\[
\Lie_X(\pi_\ast \scr{T})=\Lie_Y(\scr{T})
\]
We have canonical isomorphism of determinants:
\begin{align*}
\Delta_X^{\mathrm{add}}(\pi_\ast (-))\simeq \Delta_Y^{\mathrm{add}}(-) & \text{ on }D^b(\mathrm{CT}/Y);\\
\det_\ZZ R\Gamma_{\mathrm{W,B}}(X,\pi_\ast (-)) \simeq \det_\ZZ R\Gamma_{\mathrm{W,B}}(Y,-) & \text{ on red and blue constructible tori;}\\
\Delta_X(\pi_\ast (-))\simeq \Delta_Y(-)& \text{ on red and blue constructible tori;}\\
\end{align*}
\end{prop}
\begin{proof}
    Weil restriction is represented by the left exact functor $\pi_{\mathrm{fppf},\ast}$ which commutes with taking the Lie sheaf. Since the Lie algebra only depends on the Néron model part, we find the first equality. The second follows by taking $\ZZ$-determinants and the unique extension of determinants to $D^b(\mathrm{CT}/X)$. Since the étale realization of $\pi_\ast$ is $\pi_{\et,\ast}$, the third equality is proven in \cite{Morin2023b}. The fourth equality combines the second and third.
\end{proof}

\begin{prop}\label{delta_tensor_Q}
Let $\scr{T}^\bullet\in D^b(\mathrm{CT}/X)$. In the number field case, we have a canonical isomorphism
\[
\Delta_X^{\mathrm{add}}(\scr{T}^\bullet)_\QQ\simeq \det_\QQ \Lie_K(\scr{T}_K^\bullet)
\]
In the function field case, we have a canonical isomorphism
\[
\Delta_X^{\mathrm{add}}(\scr{T}^\bullet)_\QQ\simeq \QQ
\]
\end{prop}
\begin{proof}
    This reduces to the corresponding formula for torsion-free constructible tori, which follows in the number field case from the identification $X_\QQ=\Spec(K)$ and the commutation of the Lie algebra with flat base change, and in the function field case from the fact that coherent cohomology must be finite-dimensional over the finite base field, hence $\Lie_X(\scr{T})\otimes \QQ = 0$.
\end{proof}

\begin{prop}\label{prop:additive_part_R}
	Let $T$ be a multiplicative group over $K$ and denote $Y$ its character group. There is a canonical isomorphism
	\[
	\Lie_K(T)_\RR \simeq R\Hom_{G_\RR,X(\CC)}(Y,\CC),
	\]
	where $Y$ is seen as a $G_\RR$-equivariant sheaf on $X(\CC)$ through the decomposition groups at archimedean places of $K$.
	Let $\scr{T}^\bullet\in D^b(\mathrm{CT}/X)$, and denote $Y^\bullet$ the character group of its pullback to $K$. There is a canonical isomorphism
	\[
	\Delta^{\mathrm{add}}_X(\scr{T}^\bullet)_\RR \simeq \det_\RR R\Hom_{G_\RR,X(\CC)}(Y^\bullet,\CC).
	\]
\end{prop}
\begin{proof}
	The first isomorphism was proven in \cite[Prop. 4.8]{Morin2023b}. The second follows from the previous proposition.
\end{proof}

There is a natural $G_\RR$-equivariant map $\log|-|:\CC^\times \to \RR$ of $G_\RR$-equivariant sheaves on $X(\CC)$. Let $\scr{T}^\bullet\in D^b(\mathrm{CT}/X)$, and denote by $Y^\bullet$ the character group of its generic part (which is a bounded complex of finite type discrete $G_K$-module). Through the decomposition groups at archimedean places, $Y^\bullet$ defines a complex of $G_\RR$ equivariant sheaves on the set of complex points $X(\CC)$ equipped with the analytical topology (this is discrete and identifies with the set of complex embeddings of $K$). Denote by $\mathrm{Log}$ the composite map
\[
R\Gamma(X,\scr{T}^\bullet)\to \prod_{v\mid \infty}R\Gamma_{G_v}(K_v,\scr{T}^\bullet_v) \xrightarrow{\log} R\Hom_{G_\RR,X(\CC)}(Y^\bullet,\RR).
\]
As the target is a complex of $\RR$-vector spaces, the map factors through $R\Gamma(X,\scr{T})_{\RR}$ and we define:
\begin{defi}
    Let $\scr{T}^\bullet\in D^b(\mathrm{CT}/X)$. The Deligne compactly supported cohomology with $\RR$-coefficients of $\scr{T}^\bullet$ is defined through the fiber sequence
    \[
R\Gamma_{\mathrm{D}}(X,\scr{T}^\bullet_\RR) \to R\Gamma(X,\scr{T}^\bullet)_\RR \xrightarrow{\mathrm{Log}} R\Hom_{G_\RR,X(\CC)}(Y^\bullet,\RR).
    \]
\end{defi}
By the $\infty$-categorical nine lemma (i.e. the fiber of a morphism of fiber sequences is a fiber sequence), the fiber sequence of $G_\RR$-equivariant sheaves $(2i\pi\ZZ) \otimes \RR \to \CC \to \RR$ induces a fiber sequence
\[
R\Gamma_{\mathrm{D}}(X,\scr{T}^\bullet_\RR) \to R\Gamma_{\mathrm{B}}(X,\scr{T}^\bullet)_\RR  \to R\Hom_{G_\RR,X(\CC)}(Y^\bullet,\CC)
\]
(see \cite[§ 5]{Morin2023b}).
\begin{defi}
    Let $\scr{T}^\bullet \in D^b(\mathrm{CT}/X)$. The Weil-Arakelov cohomology with $\RR$-coefficients of $\scr{T}^\bullet$ is
    \[
    R\Gamma_{\mathrm{ar,c}}(X,\scr{T}^\bullet_\RR):=R\Gamma_{\mathrm{D}}(X,\scr{T}^\bullet_\RR)[1]\oplus R\Gamma_{\mathrm{D}}(X,\scr{T}^\bullet_\RR)
    \]
    Its $\RR$-determinant has a canonical trivialization
    \[
    \det_\RR R\Gamma_{\mathrm{ar,c}}(X,\scr{T}^\bullet_\RR) \xrightarrow{\simeq} \bigl(\det_\RR R\Gamma_{\mathrm{D}}(X,\scr{T}^\bullet_\RR)\bigr)^{-1} \otimes_\RR \det_\RR R\Gamma_{\mathrm{D}}(X,\scr{T}^\bullet_\RR) \xrightarrow{\simeq} \RR
    \]
\end{defi}

Let $\scr{T}$ be a red or blue constructible torus. As in \cite[§ 7.1]{Morin2023b}, using the duality theorem between $R\Gamma(X,\scr{T}^D)_\RR$ and $R\Gamma_{\mathrm{c,D}}(X,\scr{T}_\RR)$ given by \cite[thm. 5.2]{Morin2023b} in the number field case and \cite[thm. 4.4]{Morin2023a} in the function field case, we observe that the fiber of the composite map
\begin{align*}
R\Gamma_{\mathrm{W,B}}(X,\scr{T})_\RR &= R\Hom(R\Gamma(X,\scr{T}),\RR[-2])\oplus R\Gamma_{\mathrm{B}}(X,\scr{T})_\RR\\
& \to R\Gamma_{\mathrm{B}}(X,\scr{T})\\
& \xrightarrow{\mathrm{Log}} R\Hom_{G_\RR,X(\CC)}(Y,\CC)
\end{align*}
identifies with $R\Gamma_{ar,c}(X,\scr{T}_\RR)$, the $\RR$-determinant of which is canonically trivial. Thus using \cref{prop:additive_part_R} we obtain a natural trivialization
\[
\lambda:\Delta_X(\scr{T})_\RR\xrightarrow{\simeq} \RR.
\]
\begin{defi}
	Let $\scr{T}$ be a red or blue constructible torus. The Weil-étale Euler characteristic of $\scr{T}$ is the positive real number $\chi_X(\scr{T})$ such that
	\[
	\lambda(\Delta_X(\scr{T}))=\chi_X(\scr{T})^{-1} \ZZ \hookrightarrow \RR
	\]
\end{defi}

\begin{prop}\label{prop:functoriality_euler}
	The Weil-étale Euler characteristic is multiplicative with respect to blue-to-red short exact sequences of constructible tori. Let $Y$ be a regular integral scheme with a finite dominant morphism $\pi:Y\to X$ and let $\scr{T}$ be a red or blue constructible torus. Then
	\[
	\chi_X(\pi_\ast \scr{T})=\chi_Y(\scr{T}).
	\]
\end{prop}
\begin{proof}
	The first part is proven as in \cite[thm. 7.3]{Morin2023b} using that $\Delta^{\mathrm{add}}$ is a determinant functor and \cref{delta_tensor_Q}. The second part follows from the compatibility of the terms involved with respect to $\pi_\ast$, see \cref{compatibilty_pushforward} and \cite[prop. 7.4]{Morin2023b}.
\end{proof}

\begin{prop}
    The Weil-étale Euler characteristic has a unique extension to $D^b(\mathrm{CT}/X)$, compatible with finite dominant pushforwards.
\end{prop}
\begin{proof}
    The inclusion of red and blue constructible tori into the bounded derived category of all constructible tori induces an isomorphism on $K_0$ by the argument of \cite[def. 6.18, prop. 6.19]{Morin2023a}. The compatibility with finite dominant pushforwards then reduces to the red or blue case as in \cite[prop. 6.20]{Morin2023a}.
\end{proof}

\begin{thm}\label{thm:special_value}
	Let $X$ be a regular integral $1$-dimensional proper arithmetic scheme and let $\scr{T}^\bullet\in D^b(\mathrm{CT}/X)$. We have the vanishing order formula
	\[
	\ord_{s=0}L_X(\scr{T}^\bullet,s)=\sum (-1)^i i\cdot \dim_\RR H^i_{ar,c}(X,\scr{T}^\bullet_\RR)
	\]
	and the special value formula
	\[
	L_X^\ast(\scr{T}^\bullet,0)=\pm \chi_X(\scr{T}^\bullet).
	\]
\end{thm}

\begin{proof}
    By Artin induction and a theorem of Swan \cite[Corollary 1]{Swan1963}, this formally reduces to the case of $\GG_m$, or a torsion-free constructible torus of the form $\scr{T}=i_\ast \ZZ[-1]$ where $i:x\to X$ is a closed point. In the number field case we can then apply \cite[Thm. 7.13]{Morin2023b} to the étale realization of $\scr{T}$, whence we can conclude because the $L$-function and the Weil-étale Euler characteristic only depend on the étale realization of $\scr{T}$. In the function field case, we can proceed in those two cases by direct computation as in \cite[cor. 6.8, prop. 6.10]{Morin2023a}.
   
\end{proof}

\begin{cor}
	Let $X$ be an integral $1$-dimensional proper arithmetic scheme and let $F\in D^b(\mathbb{Z} \mhyphen \mathrm{Con} / X)$. Put $F^D:=R\sheafhom_{X_{\mathrm{et}}}(F,\GG_X)$ where $\GG_X:=[g_\ast \GG_m\to \bigoplus_{x\in X_0}i_{x,\ast}\ZZ]$ is Deninger's dualizing complex.\footnote{Such a complex would be the étale realization of an element of $D^b(\mathrm{CT}/X)$, if we had constructed such a category in the singular case.} Let $f:X\to T$ be the structural morphism to $\Spec(\ZZ)$ in the number field case or a finite morphism to some $\PP^1$ in the function field case. The complex $R\Gamma_{ar,c}(X,F^D_\RR)$ is defined as in \cite{Morin2023b}, or equivalently as $R\Gamma_{ar,c}(T,(f_\ast F)^D_\RR)$ and we have the vanishing order formula
	\[
	\ord_{s=0}L_X(F^D,s)=\sum (-1)^i i\cdot \dim_\RR H^i_{ar,c}(X,F^D_\RR)
	\]
	and the special value formula
	\[
	L_X^\ast(F^D,0)=\pm \chi_X(F^D).
	\]
	where $\chi_X:=\chi_{T}(f_\ast(-))$; equivalently, $\chi_X$ is defined by extending the above definitions for étale sheaves $F^D$ (as they only depended on the étale realization of our constructible tori), and setting $\Delta^{\mathrm{add}}_X(F^D):=\Delta^{\mathrm{add}}_{\widehat{X}}((\pi^\ast F)^D)$, where $\pi:\widehat{X}\to X$ is the normalization.
\end{cor}

\appendix

\section{$K$-theory of $\ZZ$-constructible sheaves}

\begin{prop}\label{K-theory_torsionfree}
	Let $X$ be an arithmetic scheme and let $\mathbb{Z} \mhyphen \mathrm{Con}_{\mathrm{tf}} / X\subset \mathbb{Z} \mhyphen \mathrm{Con}/ X$ denote the full subcategory of torsion-free $\ZZ$-constructible sheaves on a scheme $X$. Then the canonical inclusions induce equivalences
	\[
	K(\mathbb{Z} \mhyphen \mathrm{Con}_{\mathrm{tf}} / X) \xrightarrow{\simeq} K(\mathbb{Z} \mhyphen \mathrm{Con} / X) \xrightarrow{\simeq} K(D^b(\mathbb{Z} \mhyphen \mathrm{Con} / X))
	\]
	In particular, we have an equivalence on the $1$-truncation, which implies that those categories have the same universal determinant functor, and a determinant functor on the exact category of torsion-free $\ZZ$-constructible sheaves extends uniquely to a determinant functor on the triangulated category $D^b(\mathbb{Z} \mhyphen \mathrm{Con} / X)$.
\end{prop}
\begin{proof}
	The second equivalence is the theorem of the heart, see e.g. \cite{Barwick2015}. By \cite[\href{https://stacks.math.columbia.edu/tag/095N}{Tag 095N}]{Sta22}, any $\ZZ$-constructible sheaf receives a surjection from a torsion-free one, whose kernel must be torsion-free. The exact subcategory of torsion-free $\ZZ$-constructible sheaves is stable by subobjects and extensions so by Quillen's resolution theorem \cite[Thm. V.3.1]{Weibel2013} the left map is an equivalence.

	By the theorem of the heart for determinant functors on triangulated categories \cite{Breuning2011}, a determinant functor on $\mathbb{Z} \mhyphen \mathrm{Con}/ X$ extends uniquely to the bounded derived category; moreover by \cite{Knudsen2002}, the universal determinant functor for an exact category is computed by the $1$-truncation of $K$-theory, so a determinant functor on $\mathbb{Z} \mhyphen \mathrm{Con}_{\mathrm{tf}} / X$ extends uniquely to $\mathbb{Z} \mhyphen \mathrm{Con}/ X$.
\end{proof}

%